\documentclass[12pt,twoside,reqno]{amsart}
\usepackage[colorlinks=true,citecolor=blue]{hyperref}
\usepackage{mathptmx, amsmath, amssymb, amsfonts, amsthm, mathptmx, enumerate, color}
\usepackage[normalem]{ulem}
\usepackage{tikz,pgflibraryplotmarks}
\usepackage{fullpage}
\usepackage{graphicx}
\usepackage{mathrsfs}
\usepackage{multirow}
\usepackage{tabularx}
\usepackage{textcomp}

\usepackage{algorithm} 
\usepackage{algpseudocode}

\newcommand{\hf}{{\frac 12}}
\newcommand{\bflambda}{{\boldsymbol \lambda}}
\newcommand{\st}{\mathrm{s.t.}}

\newcommand{\bfB}{{\bf B}}

\newcommand{\bfD}{{\bf D}}

\newcommand{\bfH}{{\bf H}}
\newcommand{\bfI}{{\bf I}}

\newcommand{\bfT}{{\bf T}}

\newcommand{\bfV}{{\bf V}}
\newcommand{\bfW}{{\bf W}}
\newcommand{\bfX}{{\bf X}}

\newcommand{\bfb}{{\bf b}}

\newcommand{\bfs}{{\bf s}}
\newcommand{\bfx}{{\bf x}}
\newcommand{\bfy}{{\bf y}}
\newcommand{\bfu}{{\bf u}}

\newcommand{\bfr}{{\bf r}}

\newcommand{\bfv}{{\bf v}}

\newcommand{\edits}[1]{\textcolor{black}{#1}}

\DeclareMathOperator*{\argmin}{argmin}

\newcommand{\bbR}{{\mathbb{R}}}

\usepackage{amsmath,amsfonts,amssymb}
\usepackage{amssymb,amsfonts,amsmath,latexsym,dsfont}

\usepackage{graphicx}
\usepackage{epstopdf}

\newtheorem{theorem}{Theorem}[section]

\newtheorem{lemma}{Lemma}[section]

\theoremstyle{definition}

\numberwithin{equation}{section}

\begin{document}
\setcounter{page}{1}

\vspace*{1.0cm}
\title[Adaptive UQ-Weighted Consensus ADMM]
{Adaptive Uncertainty-Weighted ADMM for Distributed Optimization}
\author[Ye, wan, Wu Fung]{Jianping Ye$^{1,*}$, Caleb Wan$^{1,*}$, Samy Wu Fung$^2$}
\maketitle
\vspace*{-0.6cm}

\begin{center}
{\footnotesize {\it
$^1$ Department of Mathematics, University of California, Los Angeles
\\ Los Angeles, California, USA
\\
$^2$ Department of Applied Mathematics and Statistics, Colorado School of Mines
\\ Golden, Colorado, USA
}
}
\end{center}

\vskip 4mm {\small\noindent {\bf Abstract.}
    We present AUQ-ADMM, an adaptive uncertainty-weighted consensus ADMM method for solving large-scale convex optimization problems in a distributed manner. 
    Our key contribution is a novel adaptive weighting scheme that empirically increases the progress made by consensus ADMM scheme and is attractive when using a large number of subproblems. 
    The weights are related to the uncertainty associated with the solutions of each subproblem, and are efficiently computed using low-rank approximations.
    We show AUQ-ADMM provably converges and demonstrate its effectiveness on a series of machine learning applications, including elastic net regression, multinomial logistic regression, and support vector machines.
    We provide an implementation based on the PyTorch package\footnote{Code can be found at \url{http://www.github.com/chesscaleb/AUQADMM}}.
    
\noindent {\bf Keywords.}
ADMM, Preconditioning, Consensus, Machine Learning, Multinomial Logistic Regression, Support Vector Machines, Elastic Net Regression, Distributed Optimization, Convex Optimization.}

\renewcommand{\thefootnote}{}
\footnotetext{ 
$^*$ Corresponding authors. 
E-mail addresses: jianpingyemike@gmail.com (Jianping Ye), caleb@wanfamily.org (Caleb Wan).
}

\section{Introduction}
\label{sec: intro}
    We present an adaptive consensus alternating direction method of multipliers (ADMM)~\cite{boyd2011distributed} that solves 
    \begin{equation}
    	\label{eq:originalMin}
    	\begin{split}
    		\argmin_{\bfu} \;\; \sum_{j=1}^N f_j(\bfu) + g(\bfu)
    	\end{split}
    \end{equation}
    in a distributed manner, where $f_j\colon \bbR^n \mapsto \bbR$ is smooth and convex, and $g\colon \bbR^n \to \bbR$ is proximable.
    Problems of the form~\eqref{eq:originalMin} arise in many contexts, including machine learning~\cite{boyd2011distributed,xu2017adaptive2,wu2020admm,kan2020pnkh}, statistics~\cite{hastie2005elements}, phase retrieval~\cite{fung2020multigrid,candes2015phase,fienup1982phase}, geophysics~\cite{haber2014computational,fung2019multiscale}, and image processing~\cite{rudin1992nonlinear,hansen2006deblurring}. 
    These problems often contain many samples, i.e., $N$ is often very large, making the optimization computationally challenging. 
    Consensus ADMM tackles these problems by partitioning the data into $N$ smaller batches that can be solved in parallel, and in some cases, explicitly. 
    To this end,~\eqref{eq:originalMin} is reformulated as an equivalent global variable consensus problem 
    \begin{equation}
    	\begin{split}
    		\argmin_{\bfu_1,\bfu_2,\ldots,\bfu_N,\bfV} \; \quad &\sum_{j=1}^N f_j(\bfu_j) + g(\bfv)
    		\\
    		\st \quad \quad \quad &\bfu_j - \bfv = \mathbf{0},
    	\end{split}
    	\label{eq:consensusMinSum}
    \end{equation}
    where $j=1,\ldots,N$ corresponds to the different subproblems, $\bfu_j \in \bbR^n$ are the local variables, and $\bfv \in \bbR^n$ is the global variable that brings the local variables into consensus. The individual objectives $f_j$ in~\eqref{eq:consensusMinSum} are no longer coupled, allowing for the optimization problem to be solved in a distributed manner. 

    \subsection{Prior Work}
    \edits{Consensus ADMM (C-ADMM) is an ADMM-based method for solving consensus problems of the form~\eqref{eq:consensusMinSum}. A particular advantage of C-ADMM is its ease of parallelization~\cite{boyd2011distributed}, as each subproblem can be solved independently (See Sec.~\ref{sec:UQADMM}).}
    In theory, consensus ADMM converges for any positive choice of penalty parameter~\cite{eckstein1992douglas,he20121} in the augmented Lagrangian. Unfortunately, this is not always the case in practice, as the method is known to be highly sensitive to the choice of penalty parameter~\cite{ghadimi2014optimal,nishihara2015general}. To reduce the dependence of consensus ADMM on the initial penalty parameter, several adaptive penalty selection methods have been proposed. In~\cite{he2000alternating,song2016fast}, a residual balancing ADMM scheme (RB-ADMM) is proposed, which adapts the penalty parameter so that the residuals have similar magnitudes. 
    More recently, a locally adaptive consensus ADMM (AC-ADMM) scheme was proposed in~\cite{xu2017adaptive1,xu2017adaptive2}, where the penalty parameter is varied according to local curvature of the dual functions. 
    In RB-ADMM, all the local subproblems share the same penalty parameter, whereas in AC-ADMM, each subproblem acquires its own local penalty parameter. 
    Both adaptive schemes have shown to improve the performance of the standard consensus ADMM.
    
    Another drawback in consensus ADMM is that convergence of the method can deteriorate when the datasets in each batch are complementary and the number of batches $N$ is very large~\cite{fung2019uncertainty}. 
    One reason is the global averaging step performs an ``uninformed'' averaging, leading to a poor reconstruction of the global variable; see Sec.~\ref{subsec: weightIllustration}. To improve the problem of slow-averaging, a weighting scheme based on the uncertainties of the model was proposed in~\cite{fung2019uncertainty}. The weighting scheme consists of the inverse of the diagonals of the covariance matrices corresponding to each local model. In this manner, higher weights are assigned to the elements of the model in which we are more certain and vice-versa. 
    A drawback of~\cite{fung2019uncertainty} is that the weights were computed once in the off-line phase, and kept fixed throughout the optimization, i.e., the covariance matrix is assumed to be the same at each iteration.
    
    \subsection{Our Contribution}
    In this paper, we present an adaptive uncertainty-based weighting scheme that alleviates the issue of poor averaging in the global variable step in consensus ADMM. 
    We call this approach AUQ-ADMM.
    Similarly to~\cite{fung2019uncertainty}, the weights are constructed in a systematic way based on an uncertainty quantification (UQ) framework. 
    However, the AUQ-ADMM scheme presented in this paper is adaptive and extends the theory presented in~\cite{xu2017adaptive2} from scalar-weighted local models to diagonal matrix-weighted local models. 
    Following the techniques used in~\cite{xu2017adaptive2}, we provide convergence guarantees and demonstrate a convergence rate of $\mathcal{O}(1/k)$. 
    We provide an efficient GPU-based implementation in  PyTorch~\cite{NEURIPS2019_bdbca288}, a python-based library for automatic differentiation.
    
    \subsection{Outline of Paper}
    This paper is organized as follows. In Sec.~\ref{sec:UQADMM}, we introduce our proposed AUQ-ADMM method. In Sec.~\ref{sec: UQWeights}, we present a systematic way to build the weights based on uncertainties of the model. In Sec.~\ref{sec: convergence}, we show convergence guarantees of our AUQ-ADMM. In Sec.~\ref{sec: numericalResults}, we show the potential of AUQ-ADMM on a series machine learning-based tasks. We conclude with a summary and discuss future directions in Sec.~\ref{sec: conclusion}.

\section{Mathematical Derivation of AUQ-ADMM}
\label{sec:UQADMM}
    We begin by writing weighted augmented Lagrangian of~\eqref{eq:consensusMinSum}
    \begin{equation}
        \label{eq: AugmentedLagrangian}
    	\mathcal{L}(\bfv, \bfu_1, \ldots, \bfu_N, \bflambda_1, \ldots, \bflambda_N) =  g(\bfv) +
    	\sum_{j=1}^N f_j(\bfu_j) + 
    	\hf \| \bfv - \bfu_j + \bfW_j^{-1} \bflambda_j \|_{\bfW_j}^2 + 
    	\hf \| \bflambda_j \|_{\bfW_j}^2.
    \end{equation}
    where $\bfW_j \in \bbR^{n \times n}$ are positive diagonal weight matrices that spatially determine how much to constrain different areas of the local models. In the standard consensus ADMM, $\bfW_j$'s are given by the identity~\cite{boyd2011distributed}.
    

    AUQ-ADMM aims to find a saddle point of~\eqref{eq: AugmentedLagrangian} with the following iterates:
    \begin{align}
    	\bfu_j^{(k+1)} &= \argmin_\bfu \;\; \mathcal{L}(\bfv^{(k)}, \bfu_1^{(k)}, \ldots, \bfu_{j-1}^{(k)}, \bfu, \bfu_{j+1}^{(k)},\ldots, \bfu_N^{(k)}, \bflambda_1^{(k)}, \bflambda_2^{(k)}, \ldots, \bflambda_N^{(k)}), \quad \label{eq: xUpdate}
    	\\
    	&= \argmin_\bfu \;\; f_j(\bfu) + \hf \left\|(\bfv^{(k)}-\bfu) + {\bfW_j^{(k)}}^{-1} \bflambda_j^{(k)} \right\|_{\bfW_j^{(k)}}^2, \quad j=1,\ldots,N, 
    	\\
    	\bfv^{(k+1)} &= \argmin_\bfv \;\; \mathcal{L}(\bfv, \bfu_1^{(k)}, \ldots, \bfu_N^{(k)}, \bflambda_1^{(k)}, \ldots, \bflambda_N^{(k)}), \label{eq: vUpdate}
    	\\
    	&= \argmin_\bfv \;\; g(\bfv) + \hf\sum_{j=1}^N   \left\| \bfv - \bfu_j^{(k+1)} + {\bfW_j^{(k)}}^{-1} \bflambda_j^{(k)} \right\|_{\bfW_j^{(k)}}^2,
    	\\
    	\bflambda_j^{(k+1)} &= \bflambda_j^{(k)} + \bfW_j^{(k)}\left(\bfv^{(k+1)} - \bfu_j^{(k+1)}\right), \quad j=1,\ldots,N, \label{eq: lamUpdate}
    \end{align}
    where $k$ denotes the current iteration, $\bfu_j \in \bbR^n$ are the local variables, $\bfv \in \bbR^n$ is the global consensus variable, $\bflambda_j \in \bbR^n$ are the dual variables, and $\bfW_j^{(k)}$ are positive diagonal matrices (hence SPD) with norms defined as
    \begin{equation}
    	\label{eq:defWNorm}
    	\| \bfx \|_{\bfW_j^{(k)}} = \sqrt{\bfx^\top \bfW_j^{(k)} \bfx}.
    \end{equation}
    While the minimization steps in~\eqref{eq: xUpdate} are the most computationally challenging in each iteration, the costs can be alleviated by the distributed manner in which the optimization is performed.
    A further advantage is that the local subproblem can be solved using any optimization algorithm, which provides an easy way to tailor the method to different subproblems. The $\bfv$-update~\eqref{eq: vUpdate} brings the local variables $\bfu_j$ into consensus by performing a weighted averaging, and finally, the dual variables are updated.
    ascent step in~\eqref{eq: lamUpdate}.
    
    When the weight matrices $\bfW_1, \bfW_2, \ldots, \bfW_N$ are identity matrices, the iterates reduce to the standard C-ADMM algorithm, and~\eqref{eq: vUpdate} simply becomes a uniform averaging step.
    As we will show in our numerical experiments, however, when the number of splittings, $N$, is large, the performance of consensus ADMM deteriorates.  One reason is that the averaging step in~\eqref{eq: vUpdate} gives equal weighting to all elements of $\bfu_j, \;\;j=1,\ldots,N$, leading to poor reconstructions of $\bfv$, especially in the early iterations. We illustrate this in Sec.~\ref{subsec: weightIllustration}.
    
    \subsection{Stopping Criteria and Varying Penalty Parameter}
    \label{subsec: varyingPenaltyParam}
    As stopping criteria, we define the norms of the primal and dual residuals to be
    \begin{align}
    	\label{eq:normDefinitions}
    	\|\bfr^{(k)}\|_2^2 = \sum_{j=1}^N \|\bfu_j^{(k)} - \bfv^{(k)}\|_2^2, \quad \text{ and } \quad \| \bfs^{(k)} \|_2^2 = \sum_{j=1}^N \| \bfv^{(k)} - \bfv^{(k-1)} \|_{2}^2,
    \end{align}
    which are used to monitor convergence of our scheme. The iterations are terminated when 
    \begin{align}
    \label{Stopping1}
    	\|\bfr^{(k)}\|_2 \leq \epsilon_{\rm primal} \quad \text{ and } \quad \| \bfs^{(k)} \|_2 \leq \epsilon_{\rm dual},
    \end{align}
    where 
    \begin{equation}
    \label{Stopping2}
    	\begin{split}
    		\epsilon_{\rm primal} &= \sqrt{n}\epsilon_{\rm abs} + \epsilon_{\rm rel} \max\Big\{\left(\sum_{j=1}^N \|\bfu_j^{(k)}\|_2^2\right)^{1/2}, \;\; \left(N\| \bfv^{(k)} - \bfv^{(k-1)}\|_{2}^2\right)^{1/2}\Big\},
    		\\
    		\epsilon_{\rm dual} &= \sqrt{n}\epsilon_{\rm abs} + \epsilon_{\rm rel} \left(\sum_{j=1}^N \| \bflambda_j^{(k)} \|_2^2\right)^{1/2}
    	\end{split}
    \end{equation}
    are the primal and dual stopping tolerances. Here, the user must choose $\epsilon_{\rm abs}$, and $\epsilon_{\rm rel}$, which denote the absolute and relative tolerances, respectively. 
    \subsubsection{Residual Balancing (RB-ADMM)}
    \label{subsubsec: RBADMM}
    Residual balancing is a standard approach to vary the penalty parameters $\rho_j^{(k)}$ in order to improve performance~\cite{he2000alternating,boyd2011distributed}. Here, all the local subproblems share the same (constant) penalty parameter, i.e., 
    \begin{equation}
        \bfW_1^{(k)} = \bfW_2^{(k)} =\ldots=\bfW_N^{(k)} = \rho^{(k)}\bfI,
    \end{equation} 
    where $\rho^{(k)}>0$ and $\bfI \in \bbR^{n \times n}$ is the identity matrix. The scheme is defined as follows:
    \begin{align}
    	\rho^{(k)} = 
    	\begin{cases}
    		\tau \rho^{(k)} \quad &\text{ if } \| \bfr^{(k)} \|_2 > \mu \| \bfs^{(k)} \|_2
    		\\
    		\rho^{(k)}/\tau \quad &\text{ if } \| \bfs^{(k)} \|_2 > \mu \| \bfr^{(k)} \|_2
    		\\
    		\rho^{(k)} \quad &\text{ otherwise},
    	\end{cases}
    \end{align}
    where $\mu>1$ and $\tau>1$ are parameters, commonly chosen to be $\mu = 10$, and $\tau = 2$~\cite{boyd2011distributed}. The idea behind this penalty parameter update is to try to keep the primal and dual residual norms within a factor of $\mu$ of one another as they both converge to zero.
    \subsubsection{Spectral Penalty Parameter (AC-ADMM)}
    \label{subsubsec: ACADMM}
    A more recent adaptive scheme was introduced in~\cite{xu2017adaptive1} for general ADMM, and~\cite{xu2017adaptive2} for consensus ADMM (AC-ADMM). This scheme is derived from the observation that ADMM steps for the primal problem~\eqref{eq:originalMin} are equivalent to performing Douglas-Rachford splitting (DRS) on the dual formulation of~\eqref{eq:originalMin}~\cite{eckstein1992douglas}. In particular, the local penalty parameters $\rho_j^{(k)}$ are derived by assuming that the \edits{resulting subgradient function of the convex conjugates are locally linear; that is, $\partial f(x) = \alpha x + \psi$ for a set $\psi$ and scalar $\alpha$. See~\cite{xu2017adaptive2} for more details on the derivation.}
    Unlike RB-ADMM, AC-ADMM uses different parameters across local subproblems, i.e., 
    \begin{equation}
        \bfW_1^{(k)} = \rho_1^{(k)}\bfI, \quad \bfW_2^{(k)} = \rho_j^{(k+1)}\bfI, \; \ldots, \quad \bfW_N^{(k)} = \rho_N^{(k)}\bfI.
    \end{equation}
    The scheme is implemented as follows. 
    \begin{align}
    	\rho_j^{(k+1)} = \max\left( \min\Big(\hat{\rho}_j^{(k+1)}, \Big(1 + \frac{C_g}{k^2}\Big)\rho_j^{(k)} \Big), \frac{\rho_j^{(k)}}{1 + \frac{C_g}{k^2}} \right),
    \end{align}
    where 
    \begin{align}
    	\hat{\rho}_j^{(k+1)} =
    	\begin{cases}
    		\sqrt{\gamma_j^{(k)} \sigma_j^{(k)}} \quad &\text{ if } \gamma_{j,\rm{cor}}^{(k)} > \epsilon_{\rm{cor}} \;\; \text{ and } \;\; \sigma_{j,\rm{cor}}^{(k)} > \epsilon_{\rm{cor}} 
    		\\
    		\gamma_j^{(k)} \quad &\text{ if } \gamma_{j,\rm{cor}}^{(k)} > \epsilon_{\rm{cor}} \;\; \text{ and } \;\; \sigma_{j,\rm{cor}}^{(k)} \leq \epsilon_{\rm{cor}} 
    		\\
    		\beta_j^{(k)} \quad &\text{ if } \gamma_{j,\rm{cor}}^{(k)} \leq \epsilon_{\rm{cor}} \;\; \text{ and } \;\; \sigma_{j,\rm{cor}}^{(k)} > \epsilon_{\rm{cor}} 
    		\\
    		\rho_j^{(k)} \quad &\text{ otherwise }.
    	\end{cases}
    \end{align}
    Here, $\epsilon_{\rm{cor}}$ is a correlation threshold and $C_g$ is a convergence constant recommended to be chosen as $0.2$ and $10^{10}$, respectively~\cite{xu2017adaptive2}.
    Moreover, $\gamma_j^{(k)}$ and $\sigma_j^{(k)}$ are local curvature parameters, and $\gamma_{j,\rm{cor}}^{(k)}$ and $\sigma_{j,\rm{cor}}^{(k)}$ are correlation parameters which are computed using the variables at the current step $\bfu_j^{(k)}$, $\bfv^{(k)}$, $\bflambda_j^{(k)}$, and the variables at a previous step $\bfu_j^{(k_0)}$, $\bfv^{(k_0)}$, $\bflambda_j^{(k_0)}$, with $k_0$ recommended to be $k-2$. 
    For brevity and readability, we define these in the appendix.

\section{Constructing the UQ-based Weights}
    \label{sec: UQWeights}
    To represent the uncertainty of the models, the weights are chosen based on the Hessian of each individual objective.
    \begin{equation}
        \bfW_j^{(k)} \approx \sigma\left( \mathrm{diag}\left(\nabla^2 f_j\left(\bfx^{(k)}\right)\right) \right) = \sigma\left( \mathrm{diag}\left(\bfH_j^{(k)}\right) \right),
    \end{equation}
    where $\sigma$ is an adaptive affine transformation chosen to guarantee convergence (see Sec.~\ref{sec: convergence}).
    
    
    The motivation for \edits{using Hessian as the source of the weights} results from the observation that~\eqref{eq:originalMin} can be interpreted as minimizing a negative log-likehood function. 
    \edits{In the machine learning setting, $f_j$ is often the loss corresponding to a particular batch of samples, and in this case, $\bfH_j$ is the \emph{observed information matrix}~\cite{hastie2005elements,efron1978assessing}, i.e., $\bfH_j$ is an approximation of the inverse covariance matrix. See~\eqref{eq: ElasticNet},~\eqref{eq: MultiReg}, and~\eqref{eqn: SVM} for examples of $f_j$'s dependence on the data}.
    In the special case $f$ is a linear regression problem, $\bfH_j$ is exactly the inverse covariance matrix of the parameter~\cite{stuart2010inverse}.  
    
    In this UQ-framework, the diagonal elements of $\bfH_j$ represent the inverse of the variance of each element in the model. Thus, when $\bfW_j^{(k)} = \mathrm{diag}(\bfH_j^{(k)})$, \emph{ higher weights are assigned to elements in the local model with higher certainty and vice-versa}. 
    The idea for this UQ-framework was first used in the context of estimating parameters of PDEs~\cite{fung2019uncertainty}. 
    In this work, the diagonal entries of the Hessian were estimated using a low-rank approximation to alleviate computational costs.
    Moreover, the weights are built offline and \emph{kept fixed} throughout the optimization and convergence of the ADMM scheme was automatically guaranteed.
    
    Inspired by ~\cite{fung2019uncertainty}, we also construct our proposed adaptive weights using the low-rank approximation
    \begin{equation}
      \label{eq:lowRankApproxH}
      \bfW_j^{(k)} = \sigma\left( \mathrm{diag}\left(\bfV_j^{(k)} \bfD_j^{(k)} (\bfV_j^{(k)})^\top\right)\right),
    \end{equation} 
    where the right-hand-side of~\eqref{eq:lowRankApproxH} is a low-rank approximation of the Hessian with an affine transformation $\sigma$ (see Sec.~\ref{subsec: scalingDerivation}). 
    \edits{Here, $\bfV_j^{(k)} \in \bbR^{n \times r}$ is a matrix with orthonormal columns comprised of the $r$ eigenvectors corresponding the $r$ largest eigenvalues,
    $\bfD_j^{(k)} \in \bbR^{r \times r}$ is a small diagonal matrix containing the largest $r$ eigenvalues.}\\
    ~\\

    \begin{algorithm}[t]
        \caption{Restriction Interval Update}\label{alg:IntUpdt}
        \begin{algorithmic}[1]
        \State \textbf{Input: } Current iteration $k\geq1$, initial restriction interval $[a_1,b_1]$ 
        \State Set $\gamma = \frac{1}{(k+1)^2}\frac{b_1}{a_1}+1-\frac{1}{(k+1)^2}$ \Comment{Find a suitable shrinking factor $\gamma$}
        \State $a_{k+1}\gets a$ 
        \State $b_{k+1}\gets \gamma a$ 
        \State \textbf{return} $[a,b]$
        \end{algorithmic}
    \end{algorithm}
    
    \begin{algorithm}[t]
        \caption{Adaptive UQ-ADMM (AUQ-ADMM)}\label{alg:auqdmm}
        \begin{algorithmic}[1]
        \State \textbf{Input: }rank $r$ for low-rank approximation, initial $[a_1,b_1]$ and $K>0$ some integer
        \State $k=1$, $\bfv^{(0)}=\mathbf{0}$, $\bflambda_j^{(0)}=\mathbf{0}$, $\bfu_j^{(0)}$ some random vectors
        \State $a \gets a_0$, $b \gets b_0$
        \While{$k\leq\text{maxiter}$}
        \For{$j=1,2,\dots,N$}
        \State $[\bfV_j^{(k)}, \bfD_j^{(k)}] = $ lanczos($\bfH_j^{(k)}$, $r$)
        \State Update $\bfW_j^{(k)}$ using \eqref{eq:lowRankApproxH}
        \EndFor
        \State Update $[a,b]$ by \textbf{Algorithm} \ref{alg:IntUpdt}
        \For{$j=1,2,\dots,N$}
        \State Update $\bfu^{(k-1)}_j$ to $\bfu^{(k)}_j$ by (\ref{eq: xUpdate})
        \State Update $\bfv^{(k-1)}$ to $\bfv^{(k)}$ by (\ref{eq: vUpdate})
        \State Update $\bflambda^{(k-1)}_j$ to $\bflambda^{(k)}_j$ by (\ref{eq: lamUpdate})
        \EndFor
        \If {Converge by Stopping Criteria (\ref{Stopping1}) and (\ref{Stopping2})}
        \State break
        \EndIf
        \State $k\gets k+1$
        \EndWhile\label{auqadmmendwhile}
        \State \textbf{return} $\bfv^{(k)}$.
        \end{algorithmic}
    \end{algorithm}
    
    \subsection{Scaled Low-Rank Approximation}
    \label{subsec: scalingDerivation}
     The adaptive linear transformation $\sigma$ depends on a given restriction interval $[a,b]$, which restricts the range of the diagonal elements of each weight matrix $\mathbf{W}^{(k)}_j$. 
      Specifically, $\sigma(\mathbf{x})$ is defined on $[a,b]$ by the following simple element-wise linear transformation: 
    \begin{equation}
        \label{eq:adptLT}
        \begin{split}
            &\sigma(\mathbf{x}) = p\mathbf{x} + q, \quad \text{ where } \quad p=\frac{b-a}{\max{(\mathbf{x})}-\min{(\mathbf{x}})},\quad q=a-p\min{(\mathbf{x})},
        \end{split}
    \end{equation} 
    and $\max{(\mathbf{x})}$ and $\min{(\mathbf{x})}$ represent the maximum and minimum element of $\mathbf{x}$ respectively. 
    \edits{To ensure convergence of our method (see Sec.~\ref{sec: convergence}), we propose Algorithm \ref{alg:IntUpdt} to shrink the restriction interval at each iteration. The intuition is to force the weights (and hence the ADMM iterates) to converge as $k$ goes to infinity. In particular, since the diagonal entries of $\bfW_j^{(k)}$ are bounded in $[a,b]$, and $b \to a$ when $k \to \infty$ by Algorithm \ref{alg:IntUpdt}, the weights eventually converge to a constant diagonal matrix.}
    This will allow us to show asymptotic convergence of our AUQ-ADMM (see Sec.~\ref{sec: convergence})
    A key difference between our work and~\cite{fung2019uncertainty} is that our weights adaptively change at each iteration, and thus convergence is more delicate.
    
    
    \subsection{Implementation and Practical Considerations}
    \label{subsec: implementation}
    We use a Lanczos tridiagonalization method to build the eigendecomposition.
    However, we note other ways to build these low-rank approximations such as randomized methods~\cite{saibaba2016randomized} are also common. Combining the fundamentals of ADMM, UQ-based weights and restriction interval, we obtain AUQ-ADMM presented in Alg.~\ref{alg:auqdmm}.
    It is worth noting that in a parallel environment, the weights can be built locally in their corresponding processors. 
    \edits{Using the Lanczos algorithm, the constructions of the weights only require $r$ Hessian matrix-vector products, \emph{avoiding the need to build the Hessians explicitly.}
    Thus, updating $\bfu_j^{(k)}$ can be done independently and in parallel.}
    
    
    \begin{figure}[t]
    \begin{tabular}{cccc}
    $\mathbf{u}_1^{(1)}$ & $\mathbf{u}_2^{(1)}$ 
    & $\mathbf{u}_3^{(1)}$ & $\mathbf{u}_4^{(1)}$
    \\
    \includegraphics[width=0.24\textwidth]{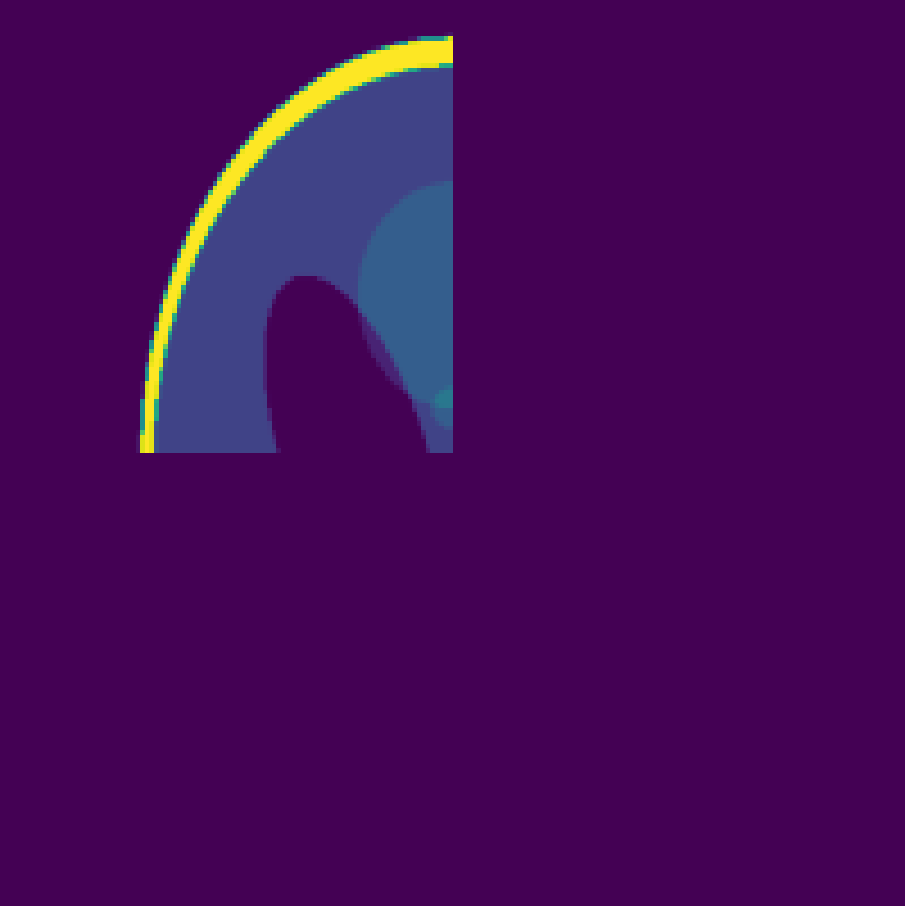}
    &
    \includegraphics[width=0.24\textwidth]{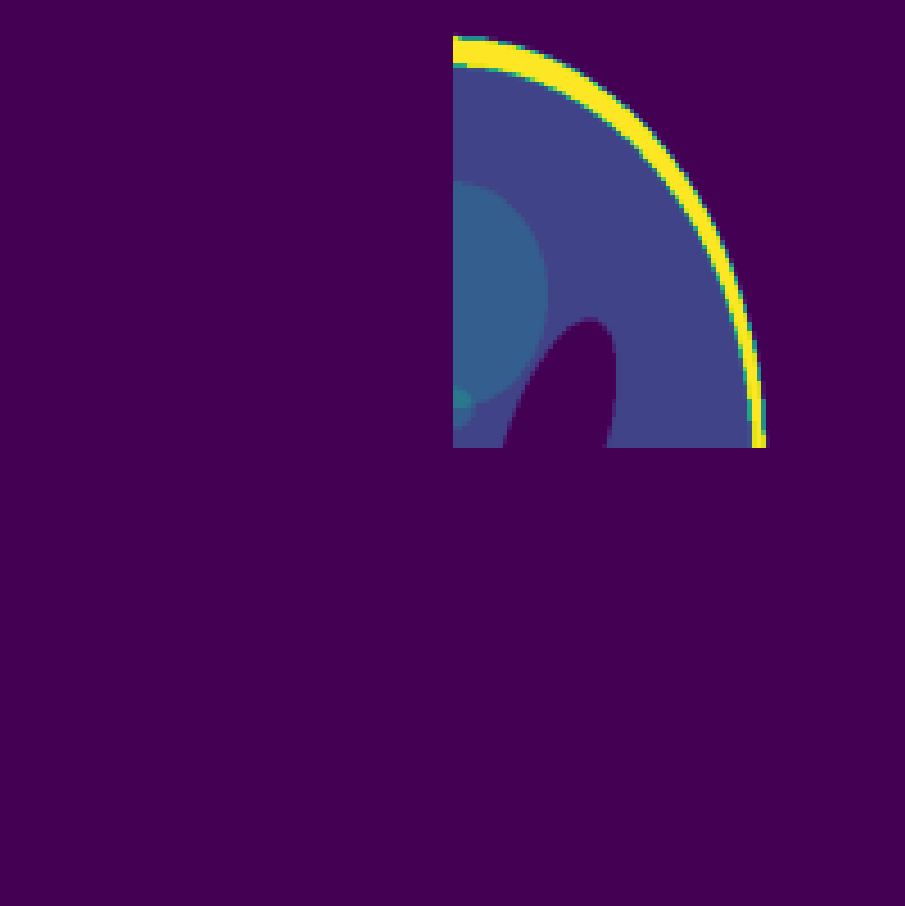}
    &
    \includegraphics[width=0.24\textwidth]{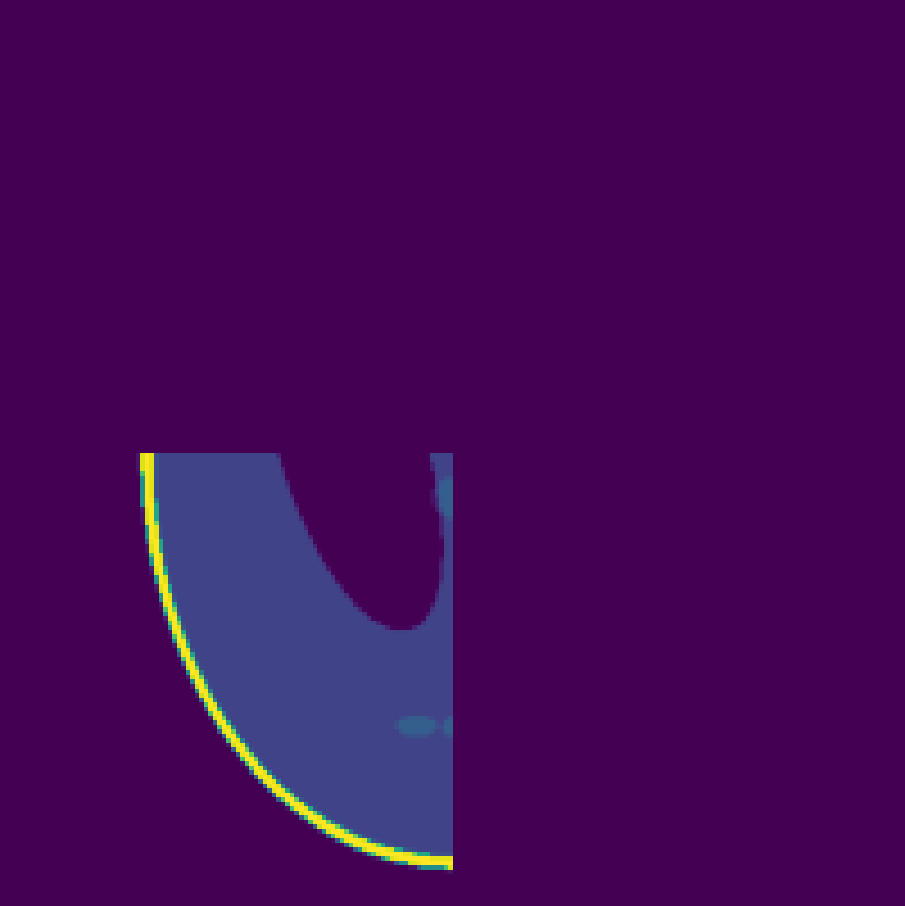}
    &
    \includegraphics[width=0.24\textwidth]{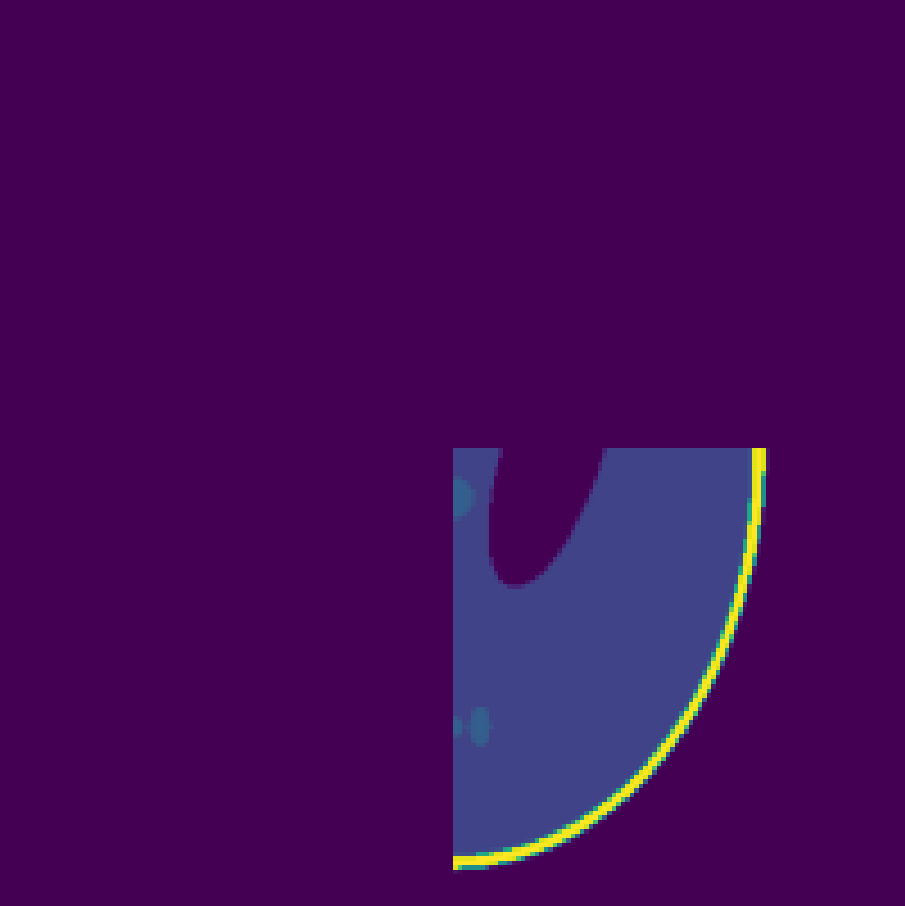}
    \\
    $\mathbf{W}_1^{(1)}$ & $\mathbf{W}_2^{(1)}$ 
    & $\mathbf{W}_3^{(1)}$ & $\mathbf{W}_4^{(1)}$
    \\
    \includegraphics[width=0.24\textwidth]{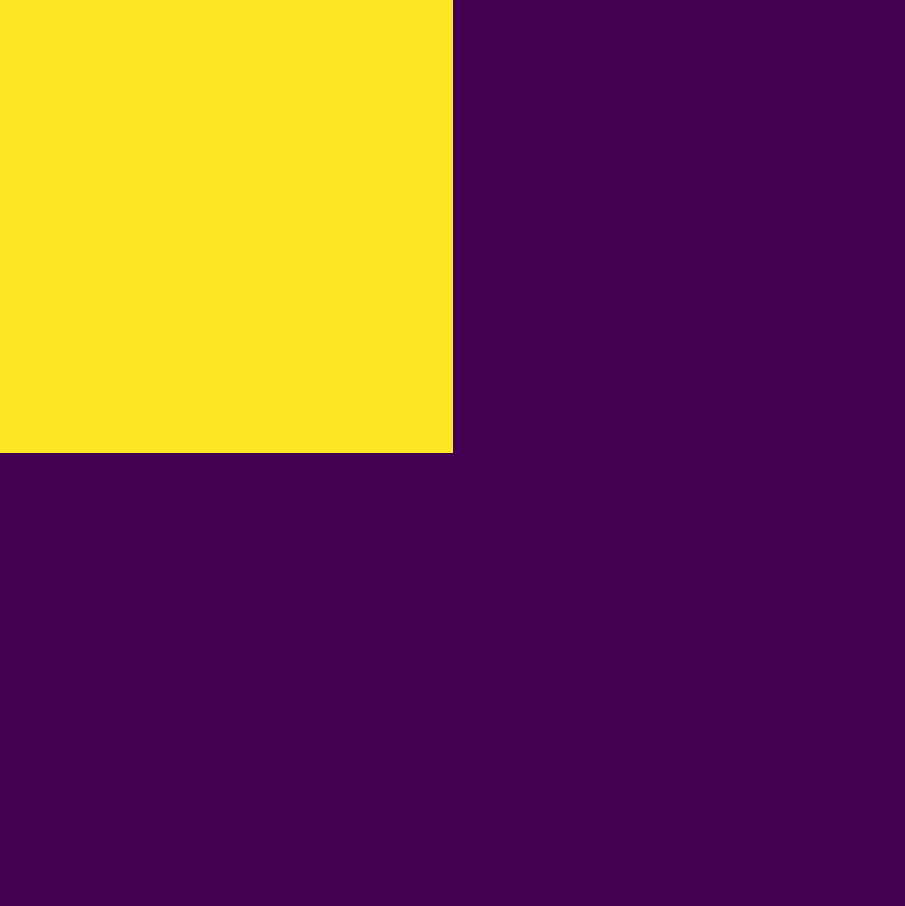}
    &
    \includegraphics[width=0.24\textwidth]{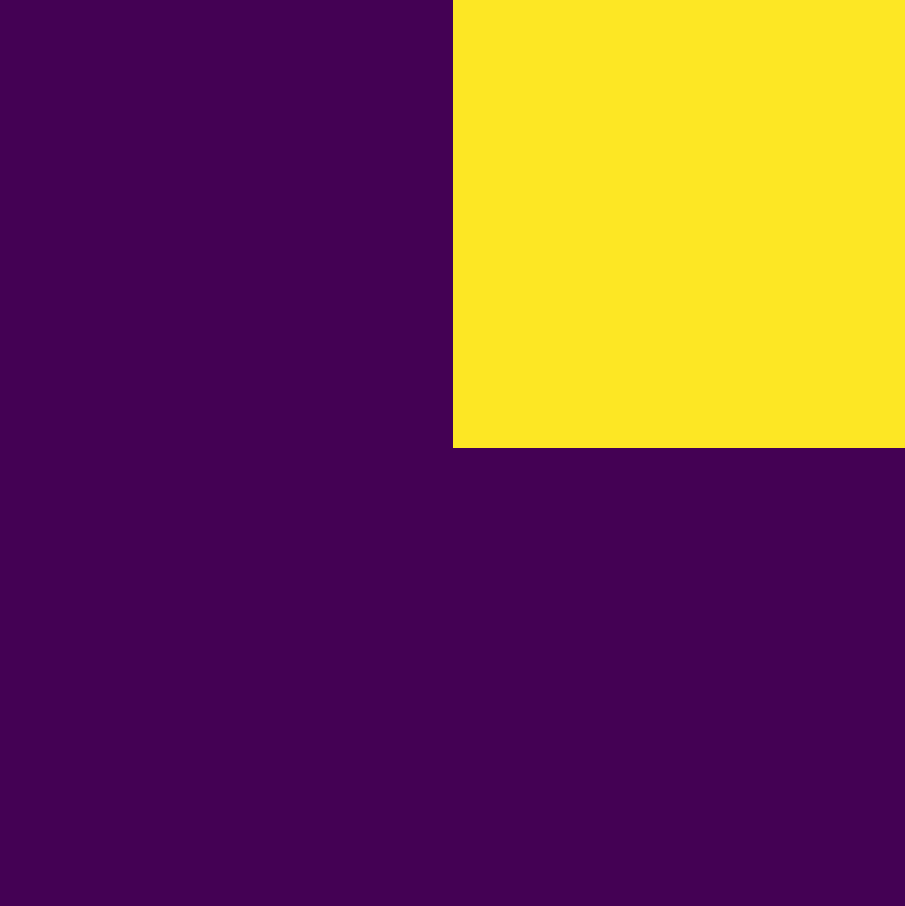}
    &
    \includegraphics[width=0.24\textwidth]{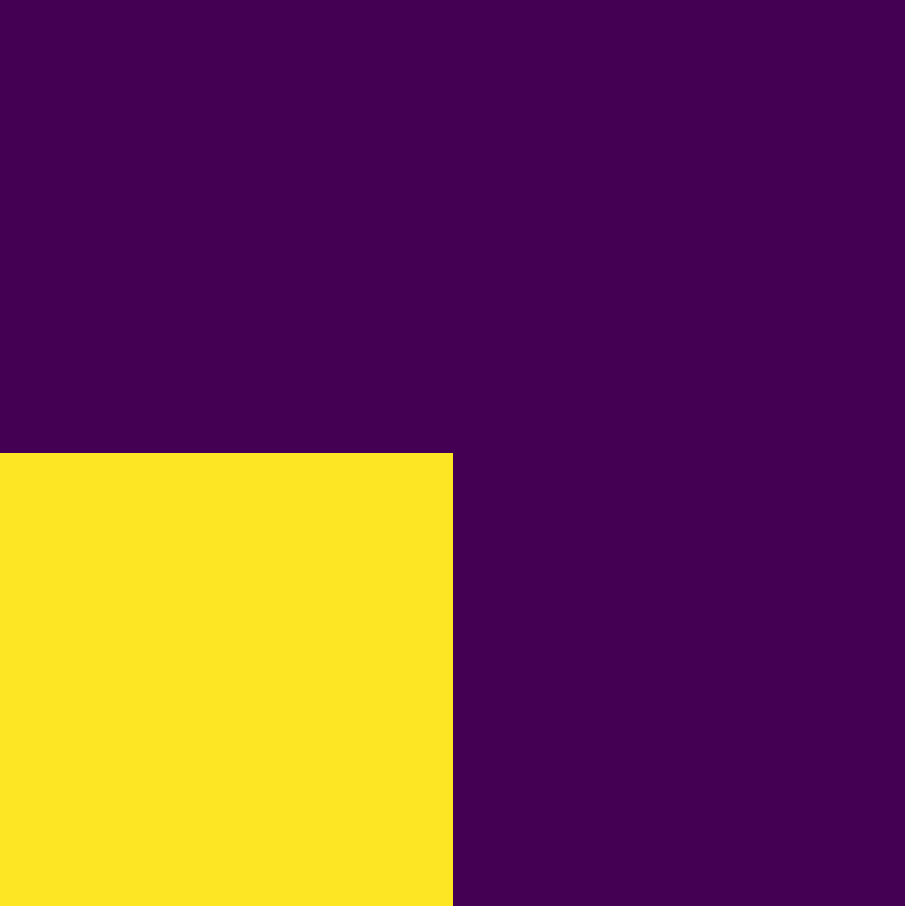}
    &
    \includegraphics[width=0.24\textwidth]{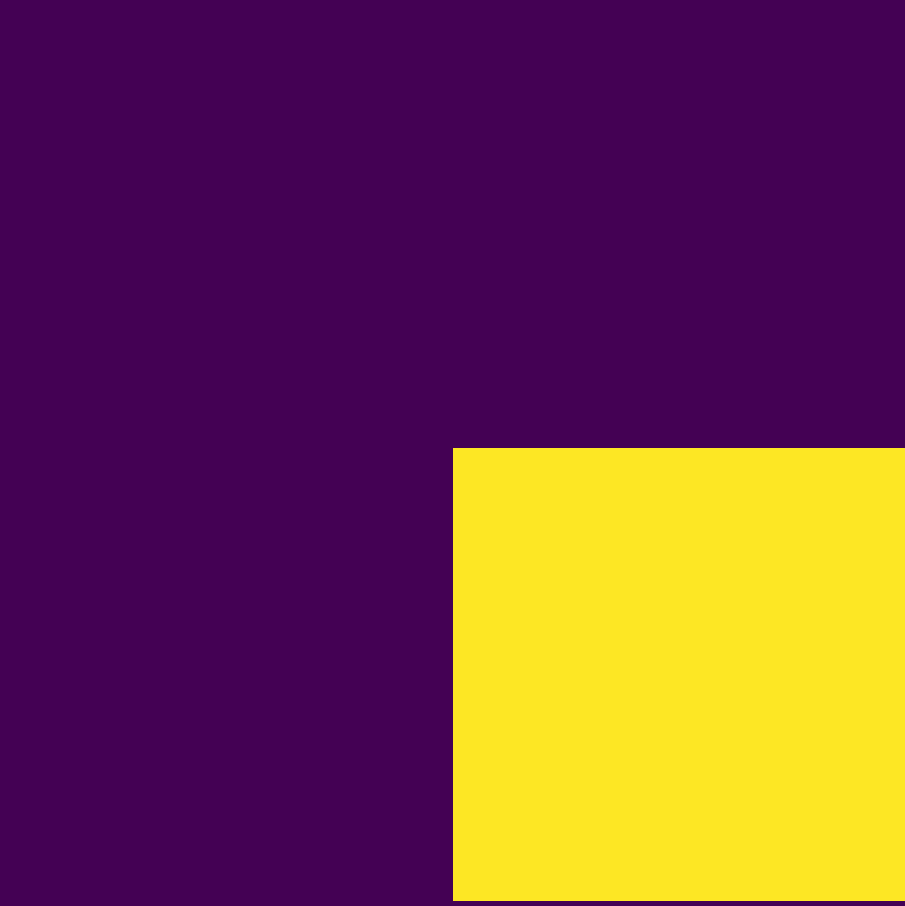}
    \\
    & $\mathbf{v}^{(1)}$ unweighted & $\mathbf{v}^{(1)}$ weighted & 
    \\
    &
    \includegraphics[width=0.24\textwidth]{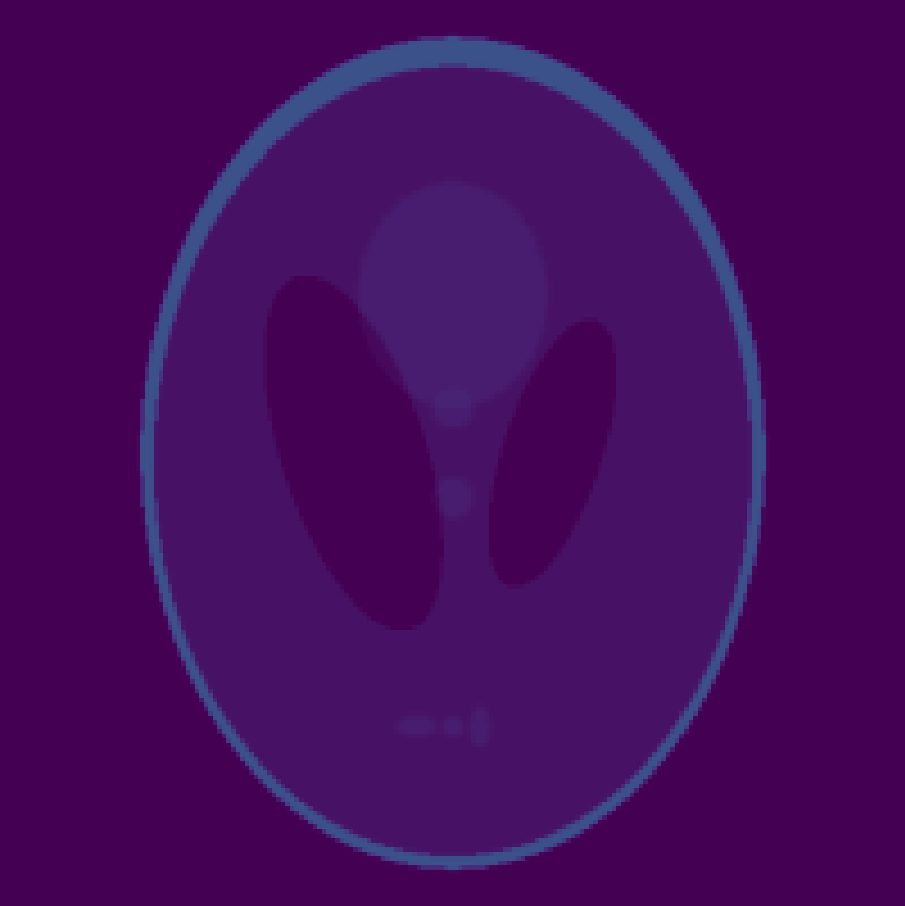}
    &
    \includegraphics[width=0.24\textwidth]{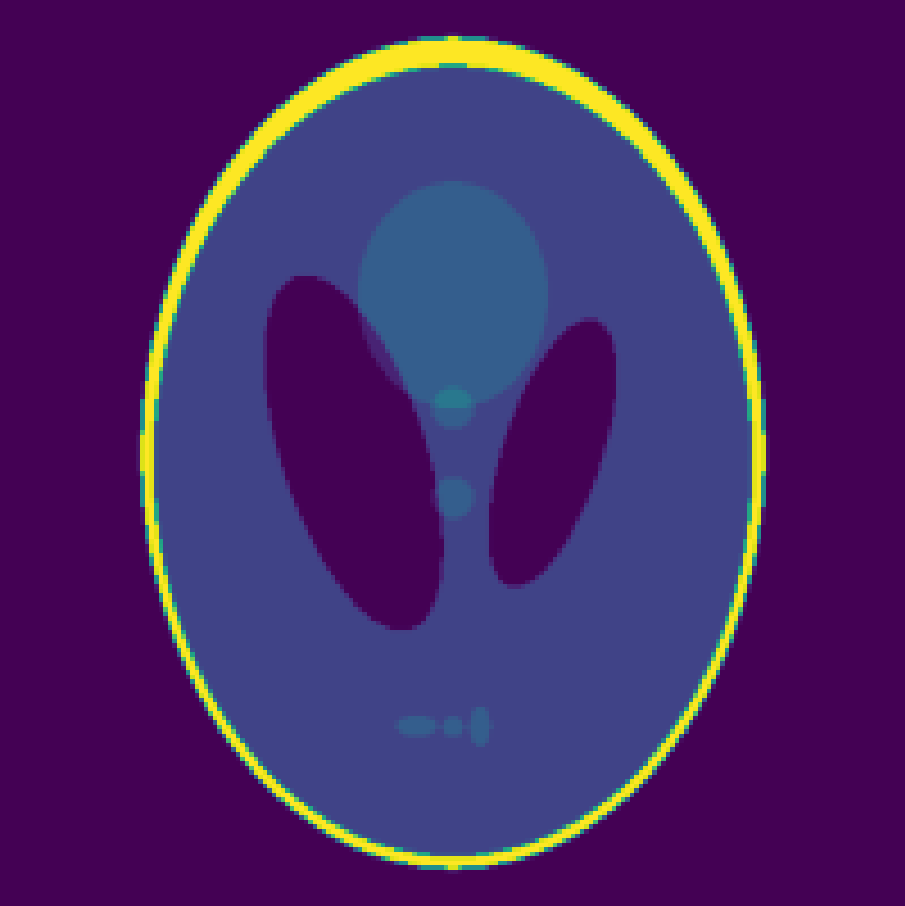}
    &
    \end{tabular}
    \caption{\small{\edits{Illustration of UQ-weighted averaging for denoising example described in Sec.~\ref{subsec: weightIllustration}. The first row shows the local models $\bfu_j$, the second row shows the \textit{diagonal} elements of the computed weights $\bfW_j$ (here, blue is a low value and yellow is a high value $~\approx 1$), \edits{and the last row show the resulting averages (i.e., $\bfv$ updates) in the first iteration. }}}}
    \label{fig:weightIllustration}
    \end{figure}
    
    
    \subsection{Illustration of Weighted Averaging}
    \label{subsec: weightIllustration}
    To illustrate the intuition behind the weighting scheme, we consider a 2D image denoising example from given by
    \begin{equation}
      \argmin_{\bfu} \;\; \hf \| \bfu - \bfb \|_2^2 + \frac{\alpha}{2}\|\bfu\|_2^2,
    \end{equation}
    where $\bfb \in \bbR^n$ is the observed noisy image, and $\alpha = 10^{-3}$. 
    The corresponding consensus problem is given by 
    \begin{equation}
      \begin{split}
        \argmin_{\bfu_1,\ldots,\bfu_N,\bfv} \quad &\sum_{j=1}^N \hf \| \bfu_j - \bfb_j \|_2^2 + \frac{\alpha}{2}\|\bfv\|_2^2
        \\
        \st \quad \quad & \bfu_j - \bfv = \mathbf{0},
      \end{split}
    \end{equation}
    where we use $N=4$ splittings corresponding to the $4$ quadrants of the image. 
    Here, \edits{$\bfb_j \in \bbR^{n}$ corresponds to modified versions of the original image containing the same pixel values in the $j^{th}$ quadrant and zeroes in the remaining quadrant (see Fig.~\ref{fig:weightIllustration}).}
    
    Fig.~\ref{fig:weightIllustration} compares the averaged reconstructions in the first iteration of the unweighted and UQ-weighted scheme. 
    In particular, Fig.~\ref{fig:weightIllustration} shows that introducing the weights considerably improves the quality of the averaging step in consensus ADMM. 
    \edits{This is because standard C-ADMM uniformly averages all $\bfu_j$ whereas AUQ-ADMM uses the uncertainty-based weights $\bfW_j$ to perform a weighted averaging.
    We note that even though the weights could be intuitively constructed by hand for this particular example; in general, it is not clear how to manually design the weights.
    This example shows that we have a principled way to construct the weighting scheme that agrees with intuitively/manually choosing the weights in obvious cases.}
    
\section{Convergence}
\label{sec: convergence}
Similar to~\cite{xu2017adaptive2}, we begin with some useful notation. Let 
\begin{align}
    \bfB&=-(\bfI_n;\dots;\bfI_n)\in\mathbb{R}^{Nn\times n}, \\
    \bfy &= (\bfu;\bfv;\bflambda) \in \bbR^{(2N+1)n}, \\
  \bfW^{(k)} &= 
  \begin{bmatrix} 
    \bfW^{(k)}_{1} 
    \\
    & \bfW^{(k)}_{2}
    \\
    & & \ddots
    \\
    & & & \bfW^{(k)}_{N}
  \end{bmatrix}
  \in \bbR^{Nn \times Nn}.
\end{align}
Denote the ADMM iterates by
\begin{align}
    \bfy^{(k)}&=(\bfu^{(k)};\bfv^{(k)};\bflambda^{(k)}),
\end{align}
and let
\begin{align}
    \hat{\bflambda}^{(k)}&=\bflambda^{(k-1)}+\bfW^{(k-1)}(-\bfu^{(k)}-\bfB\bfv^{(k-1)}),
    \\
    \Tilde{\bfy}^{(k)}&=(\bfu^{(k)};\bfv^{(k)};\hat{\bflambda}^{(k)}). 
\end{align}
Finally, set 
\begin{align}
    \phi(\bfu,\bfv) = f(\bfu)+g(\bfv) \notag, \quad F(\bfy) = (-\bflambda;-\bfB^T\bflambda;\bfu+\bfB\bfv), \notag
\end{align}
and 
\begin{equation}
  \bfT^{(k)} = 
  \begin{pmatrix} 
    \mathbf{0}
    \\
    & \bfB^T\bfW^{(k)}\bfB
    \\
    & & (\bfW^{(k)})^{-1}
  \end{pmatrix}\in\bbR^{(2N+1)n \times (2N+1)n}. \notag
\end{equation}
\begin{lemma}
    \label{lemma: WeightConvg}
        In Algorithm ~\ref{alg:IntUpdt} and Algorithm~\ref{alg:auqdmm}, the weights $\{\bfW^{(k)}\}_{k=0}^{\infty}$ ($\bfW^{(0)}=\bfW^{(1)}$) generated by AUQ-ADMM satisfy
        \begin{align}
            \bfW^{(k)} \preceq (1+c^{(k)}) \bfW^{(k-1)} \; \text{ and } \; (\bfW^{(k)})^{-1} \preceq (1+c^{(k)}) (\bfW^{(k-1)})^{-1}, \label{tempEQ}
          \end{align}
        where
        \begin{equation}
        \label{eqn: ck}
        c^{(k)} = 
        \begin{cases}
          \frac{b_1}{a_1}-1& \text{if $k=1$}
          \\
          \frac{b_{k-1}}{a_{k-1}}-1. & \text{if $k>1$}
        \end{cases},
        \qquad k = 1,2,3,\ldots
        \end{equation}
        is a sequence of positive scalars satisfying
        \begin{equation}
            \sum_{k=1}^\infty c^{(k)} < \infty \notag,
        \end{equation}
        and $[a_1, b_1]$ is the initial restriction interval. Here, the matrix inequality $\preceq$ in~\eqref{tempEQ} implies 
          \begin{equation}
            \forall \bfx, \; \| \bfx \|_{\bfW^{(k)}}^2 \leq (1+c^{(k)}) \| \bfx \|_{\bfW^{(k-1)}}^2, \; \text{ and } \; \| \bfx \|_{(\bfW^{(k)})^{-1}}^2 \leq (1+c^{(k)}) \| \bfx \|_{(\bfW^{(k-1)})^{-1}}^2. \notag
        \end{equation}
    \end{lemma}
    
    
    \begin{theorem}
        \label{theorem: Convg}
        In the AUQ-ADMM scheme (Algorithm \ref{alg:auqdmm}), the sequence $\bar{\bfy}^i$ defined as $\bar{\bfy}^i=\frac{1}{i}\sum_{k=1}^{i}\Tilde{\bfy}^{(k)}$ generated by the AUQ-ADMM satisfies
        \begin{equation}
          \phi \Big( \frac{1}{i}\sum_{k=1}^{i}\bfu^{(k)},\frac{1}{i}\sum_{k=1}^{i}\bfv^{(k)} 
          \Big)
          -\phi(\bfu^*,\bfv^*)+(\bar{\bfy}^i-\bfy^*)F(\bar{\bfy}^i)\leq\frac{\|\bfy^*-\bfy^{(0)}\|^2_{\bfT^{(0)}}+C\|\bfy^*-\bfy^{(1)}\|^2_{\bfT^{(0)}}}{2i},
        \end{equation}
        where $\bfy^*=(\bfu^*;\bfv^*;\bflambda^*)$ is optimal, and
        $$C=\left(\sum_{k=1}^\infty c^{(k)}\right) \left(\prod_{k=1}^{\infty}(1+c^{(k)})\right)<\infty,$$
        where $c^{(k)}$ is defined by (\ref{eqn: ck}).
    \end{theorem}
    
    Lemma \ref{lemma: WeightConvg} shows the convergence of weights $\bfW^{(k)}$, and Theorem \ref{theorem: Convg} shows the $\mathcal{O}(1/k)$ ergodic convergence rate. Proofs can be found in in App. \ref{Appendix: A}.

\section{Numerical Results}
    \label{sec: numericalResults}
    In this section, we outline the potential of our proposed scheme. 
    We test our proposed AUQ-ADMM on a series of machine learning tasks, including elastic net regression~\cite{hastie2005elements}, multinomial logistic regression~\cite{hastie2005elements}, and support vector machines (SVMs)~\cite{cortes1995support}.
    
    \subsection{Experimental Setup} 
    \label{subsec: setup}
    In our experiments, we compare five algorithms: our proposed adaptive uncertainty-weighted consensus ADMM (AUQ-ADMM), consensus ADMM (C-ADMM)~\cite{boyd2011distributed}, residual-based ADMM (RB-ADMM)~\cite{boyd2011distributed}, and adaptive consensus ADMM (AC-ADMM)~\cite{xu2017adaptive2}.
    For all algorithms, we use initial penalty parameter $\tau_0 = 1$, absolute stopping tolerance $\epsilon_\text{abs}=10^{-4}$, relative stopping tolerance  $\epsilon_\text{rel}=10^{-5}$, and a maximum of 250 iterations. 
    
    \subsection{Datasets}
    \label{subsec: datasets}
\begin{figure}[htb]
      \label{fig:datasets}
      \centering
      \begin{tabular}{ccc}
      MNIST & CIFAR10 & SVHN
      \\
      \includegraphics[width=0.28\textwidth, height=1.6in]{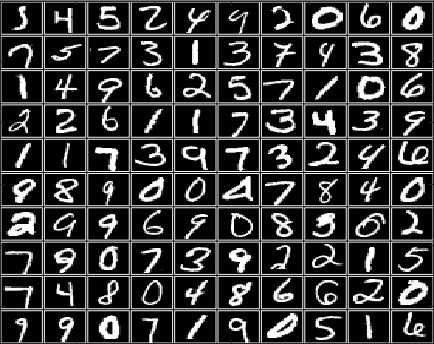}
      &
      \includegraphics[width=0.28\textwidth, height=1.6in]{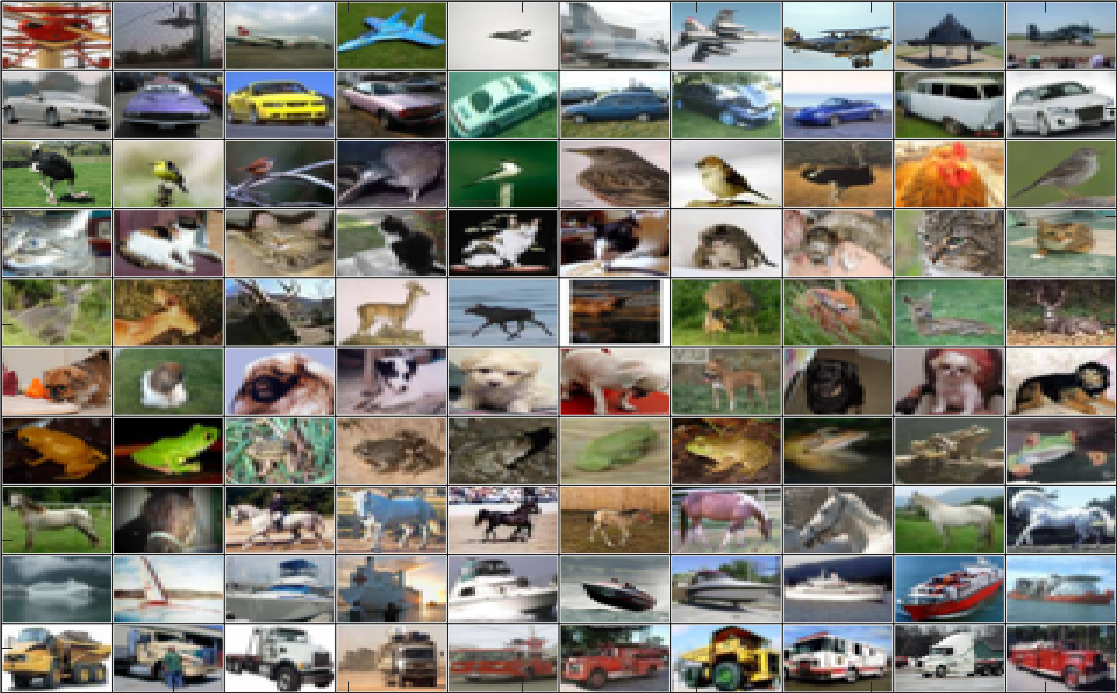}
      &
      \includegraphics[width=0.28\textwidth,
      height=1.6in]{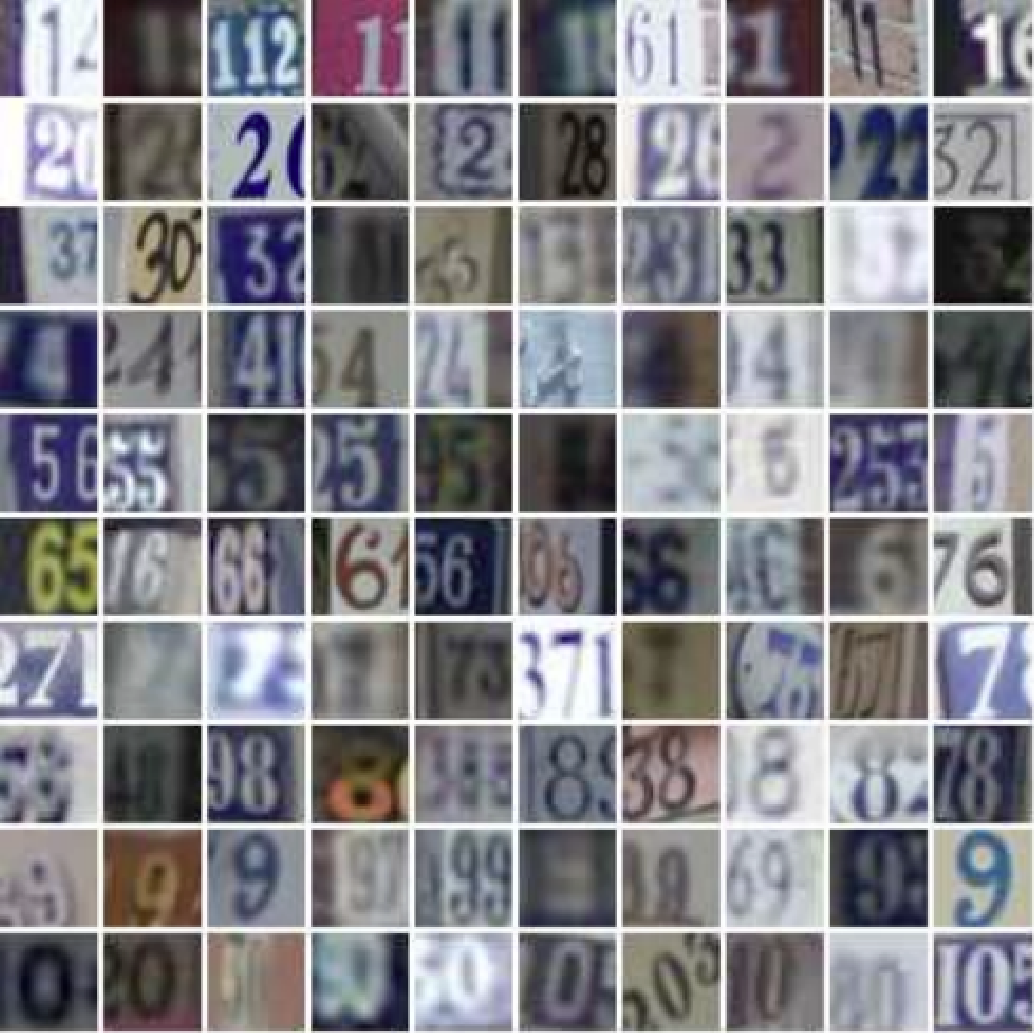}
      \end{tabular}
      \caption{\small{Sample images from MNIST dataset (left) and CIFAR10 dataset (center), and SVHN (right) described in Sec.~\ref{subsec: datasets}}}
    \end{figure}

    We also use three benchmark datasets for image classification: MNIST~\cite{lecun1998mnist}, CIFAR10~\cite{krizhevsky2009learning}, and SVHN~\cite{netzer2011reading}. 
    The MNIST dataset consists of $60,000$ labeled digital images of size $28 \times 28$ showing hand written digits from $0$ to $9$. 
    The CIFAR10 dataset consists of $60,000$ RGB images of size $32 \times 32$ that are divided into $10$ classes: airplane, automobile, bird, cat, deer, dog, frog, horse, ship, and truck.
    Finally, the SVHN dataset consists of $60,000$ RGB images obtained from house numbers in Google Street View images.
    Examples of these images are shown in Fig.~\ref{fig:datasets}.
    
    \subsection{Varying the Rank}
    We compare convergence results of AUQ-ADMM on various ranks used in the low rank approximation step. 
    Here, compute the weights for multinomial regression using hte MNIST dataset. 
    The restriction interval is initialized to be $[0.1,1.0]$, 2000 samples per worker, one class per worker, 10 workers in total. 
    \begin{figure}[htb]
		\centering
		\setlength{\tabcolsep}{1pt}
		\small
		\begin{tabular}{ccc}
		    Loss $f+g$ & Primal Residual & Dual Residual
			\\
			\includegraphics[width=0.3\textwidth]{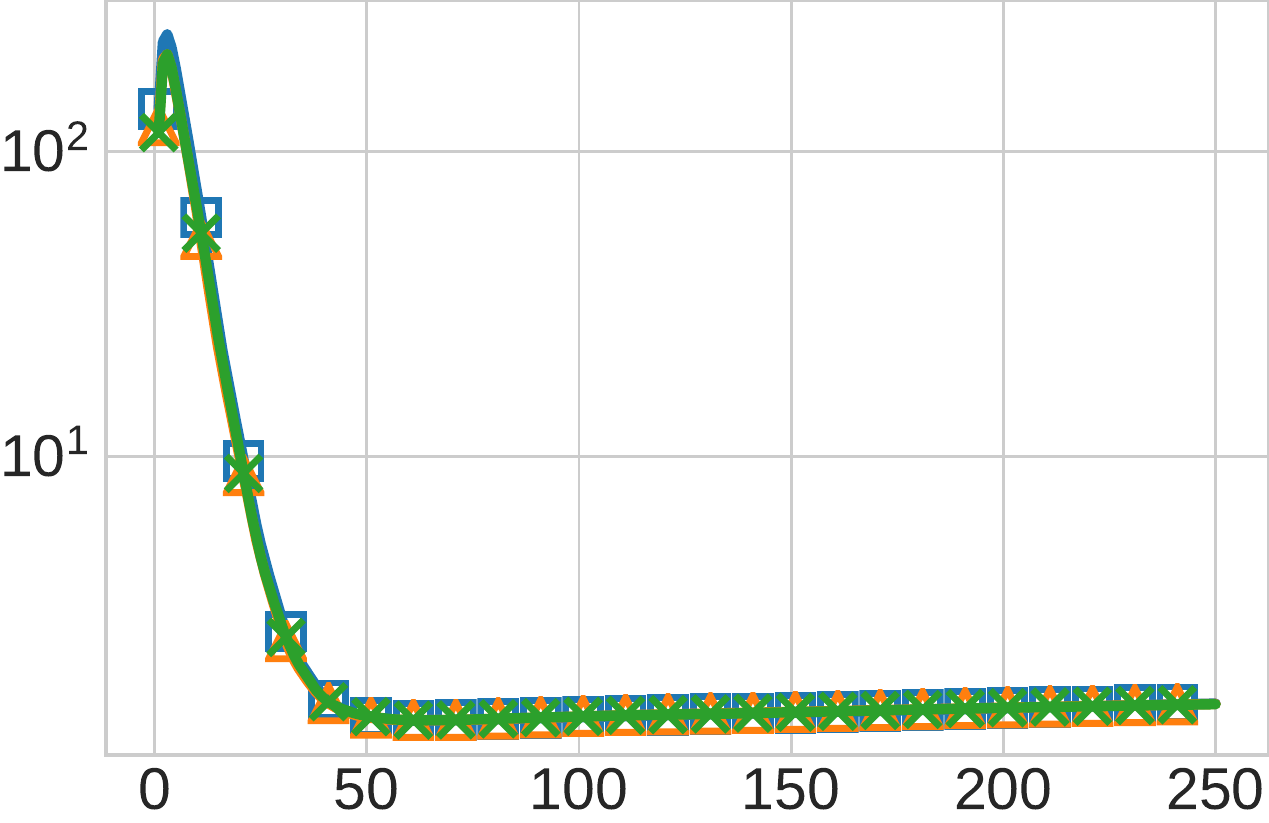}
			&
			\includegraphics[width=0.3\textwidth]{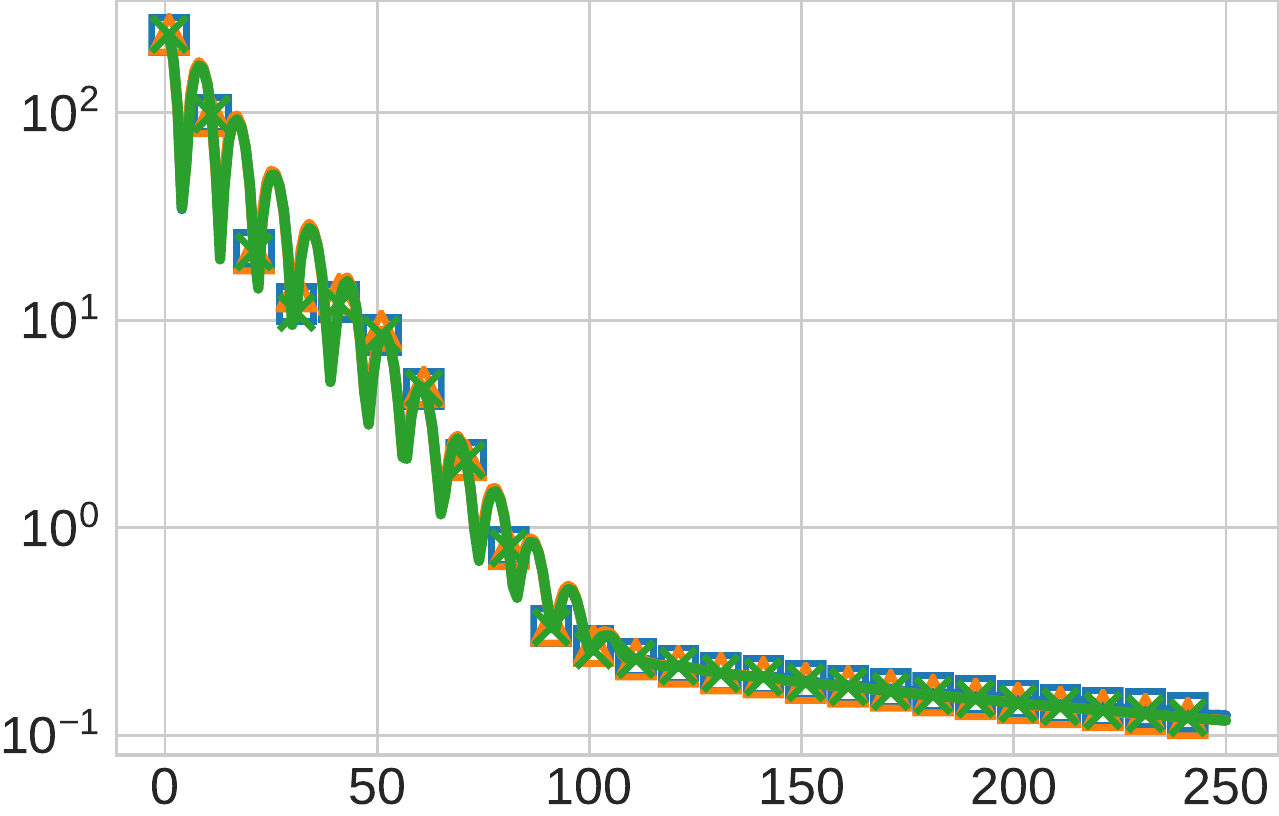} 
			&
			\includegraphics[width=0.3\textwidth]{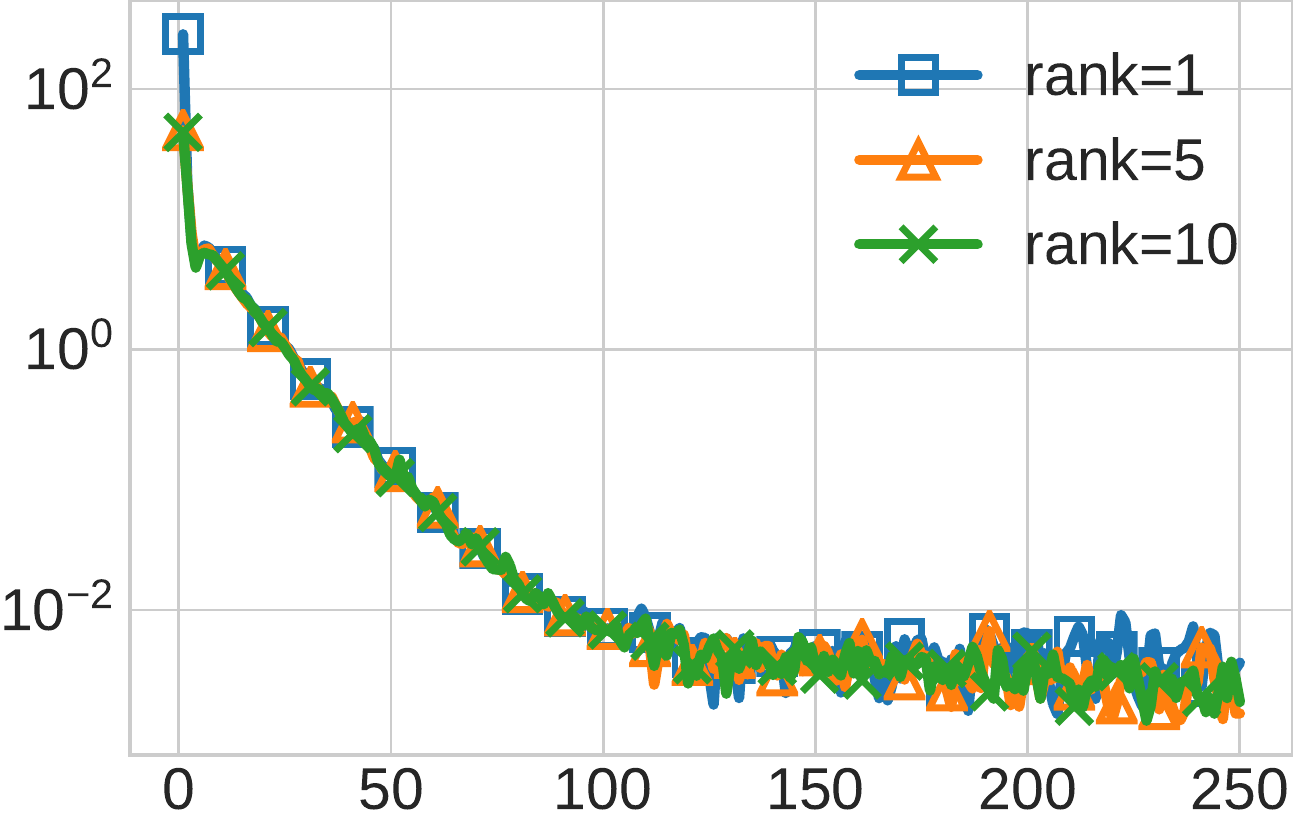}
			\\
			iterations & iterations & iterations 
		\end{tabular}
		\caption{Convergence of AUQ-ADMM under different ranks, with multinomial logistic regression on MNIST dataset. Ranks chosen are: 1, 5, 10 and 100.}
		\label{fig:Rank}
	\end{figure}
	    
	As shown in Figure \ref{fig:Rank}, there are no differences in the results produced by the different ranks. 
    This suggests that the data has an inherently low-rank structure, and allows us to compute the weights at very low computational costs.
	Motivated by this observation, we use a rank 5 approximation of the Hessian in the remainder of our experiments.
    
    \subsection{Elastic Net Regression}
    \label{subsec: elastic_net_regression}
\begin{figure}[htb]
		\centering
		\setlength{\tabcolsep}{1pt}
		\small
		\begin{tabular}{ccc}
		    MNIST Loss $f+g$ & SVHN Loss $f+g$ & CIFAR10 Loss $f+g$
			\\
			\includegraphics[width=0.3\textwidth]{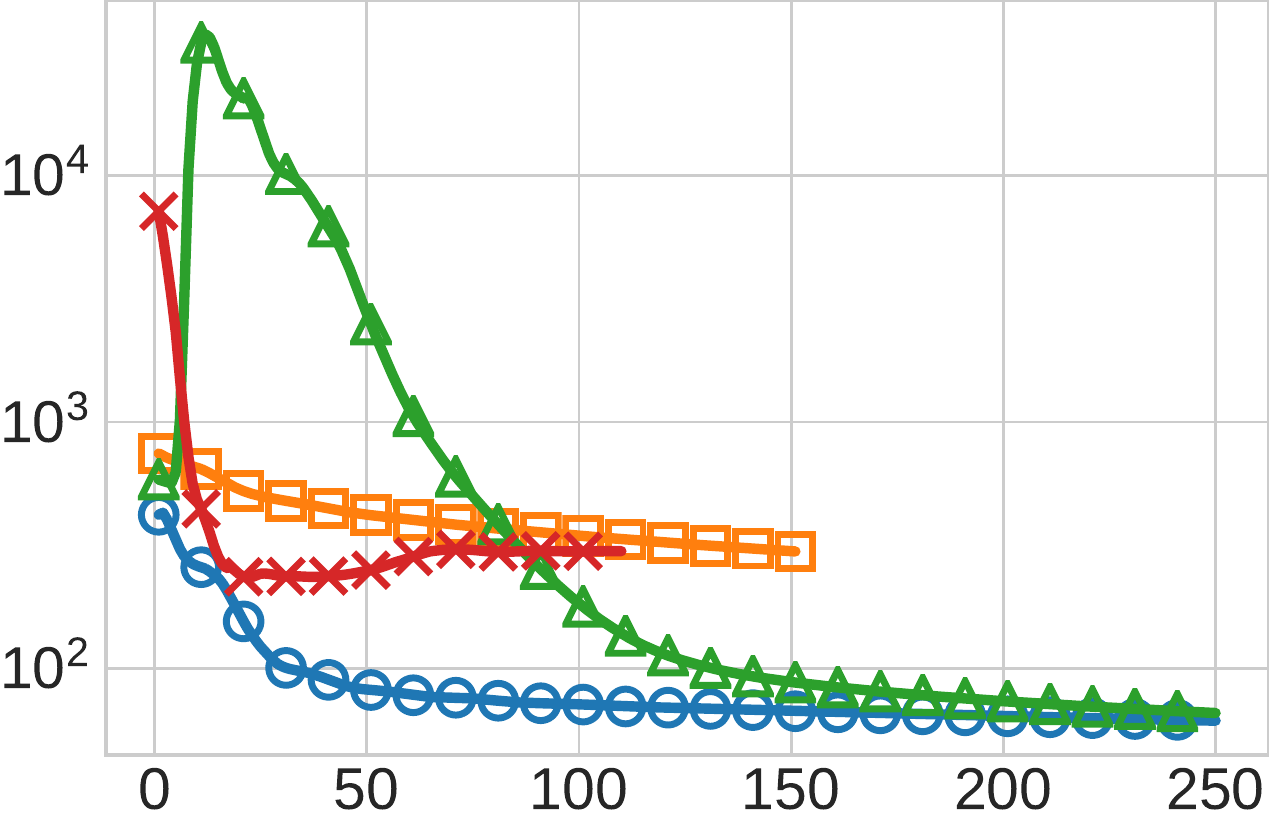}
			&
			\includegraphics[width=0.3\textwidth]{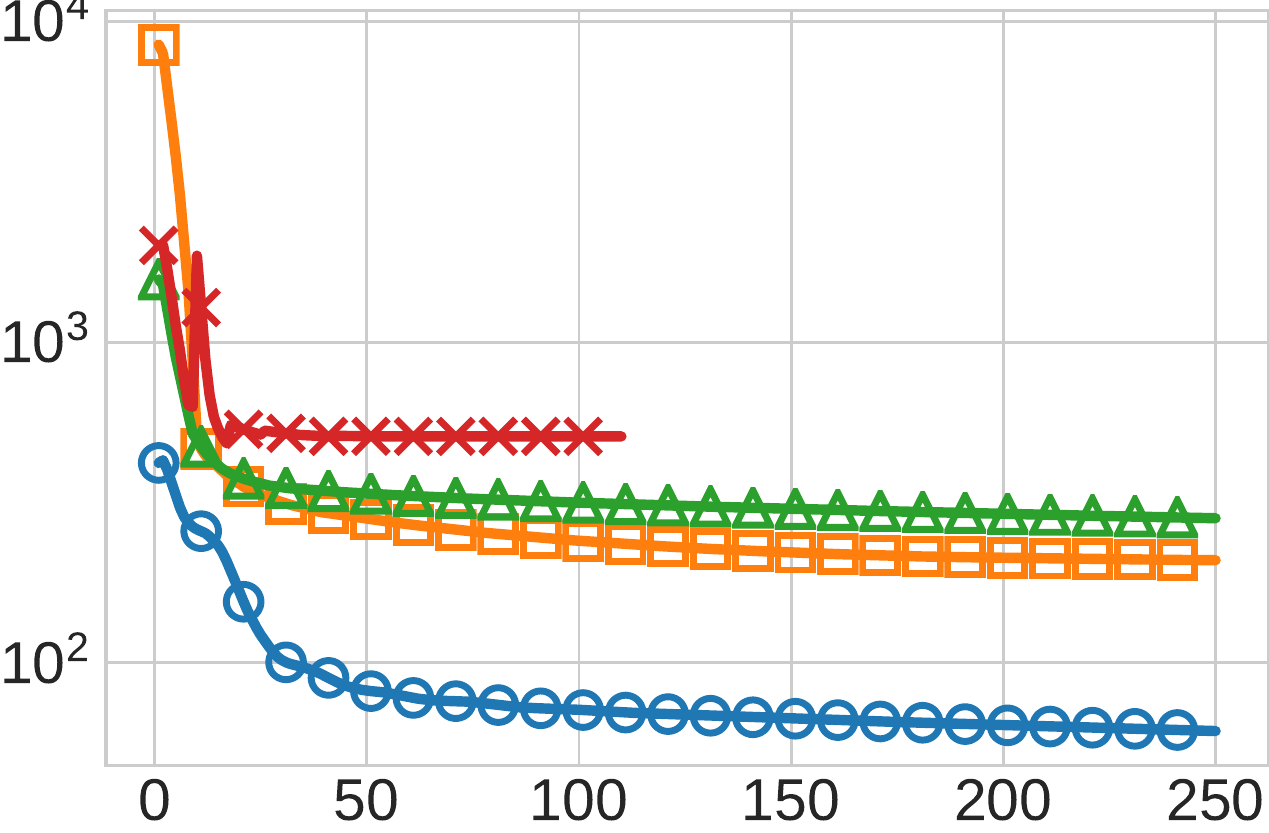} 
			&
			\includegraphics[width=0.3\textwidth]{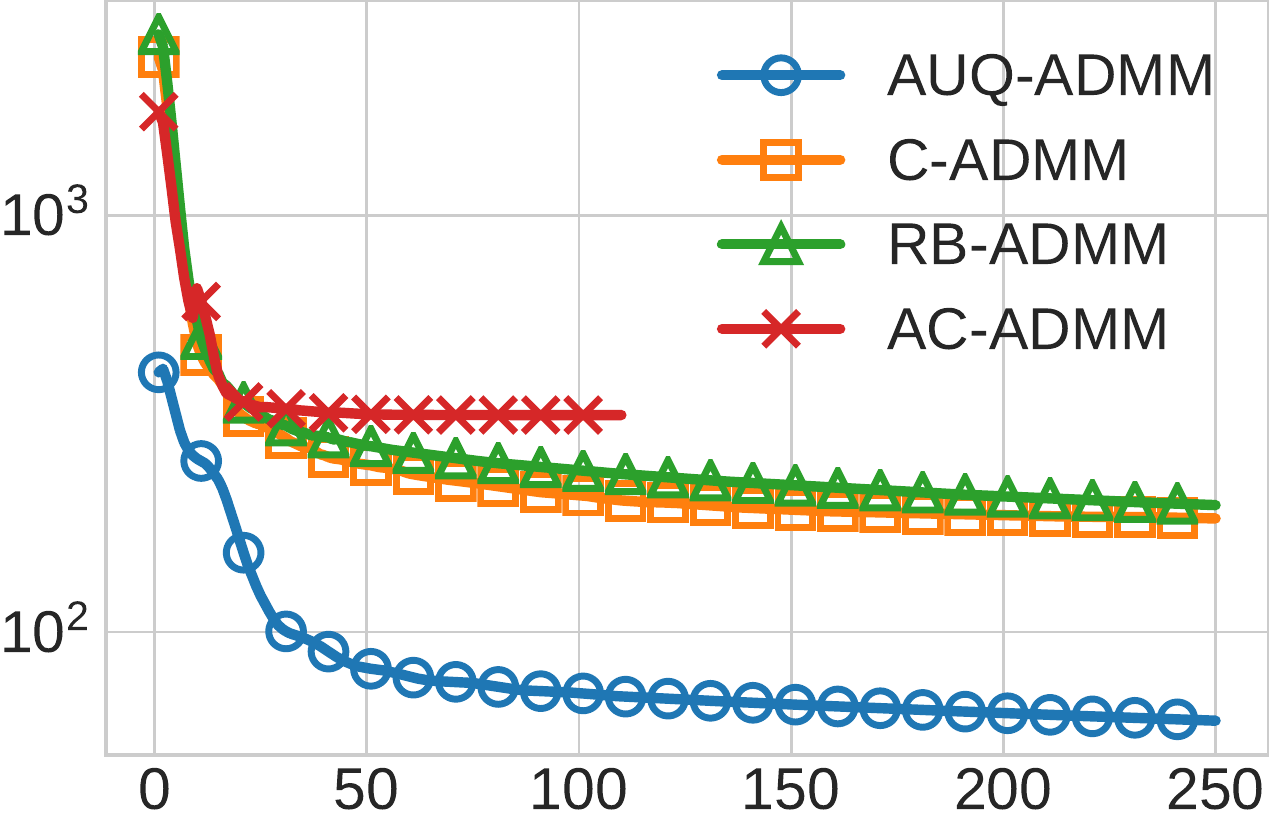}
			\\
			iterations & iterations & iterations
		\end{tabular}
		\caption{
        Linear Regression with Elastic Net Regularizer with MNIST and CIFAR10: loss function value comparison among different ADMM-algorithms.
		}
		\label{fig:Elastic}
	\end{figure}

    We consider the linear regression with the elastic net regularizer. More specifically, $f_j(\bfu_j)$ and $g(\bfv)$ are defined as: 
    \begin{equation}
        \label{eq: ElasticNet}
        f_j(\bfu_j) = \frac{1}{2}||\mathbf{X}_j\mathbf{u}_j-\bfy_j||^2,\quad g(\bfv)={\rho_1}|\bfv| + \frac{\rho_2}{2}\|\bfv\|^2_2,
    \end{equation}
    where $\bfu_j,\bfv\in\bbR^{m\times 1}$, $\mathbf{X}_j\in\bbR^{n_j\times m}$ ($n_j$ is the sample size for worker $j$) is the input data matrix for worker $j$, $\bfy_j\in\bbR^{n_j}$ is the corresponding labels and $\rho_1 = \rho_2 = 10^{-2}$ are regularization parameters.
    In our setup, $j\in\{1,2,\dots,10\}$; $m=784$, $n_j=2000$ for MNIST; $m=3072$, $n_j=1000$ for CIFAR10 and SVHN, and the data are specifically divided such that each worker only processes one class. In other words, MNIST has 2000 samples per worker, one class per worker; SVHN and CIFAR10 have 1000 samples per worker, one class per worker, 10 workers in total. The restriction interval is initialized to be $[0.1,1]$. The convergence results are shown in Figure \ref{fig:Elastic}.
	
	In all three cases, the proposed AUQ-ADMM has the potential to outperform all other ADMM algorithms performs competitively. 
    
    \subsection{Multinomial Logistic Regression}
    \label{subsec: multinomial_logistic_regression}
\begin{figure}[htb]
		\centering
		\setlength{\tabcolsep}{1pt}
		\small
		\begin{tabular}{ccc}
		    MNIST Loss $f+g$ & SVHN Loss $f+g$ & CIFAR10 Loss $f+g$
			\\
			\includegraphics[width=0.3\textwidth]{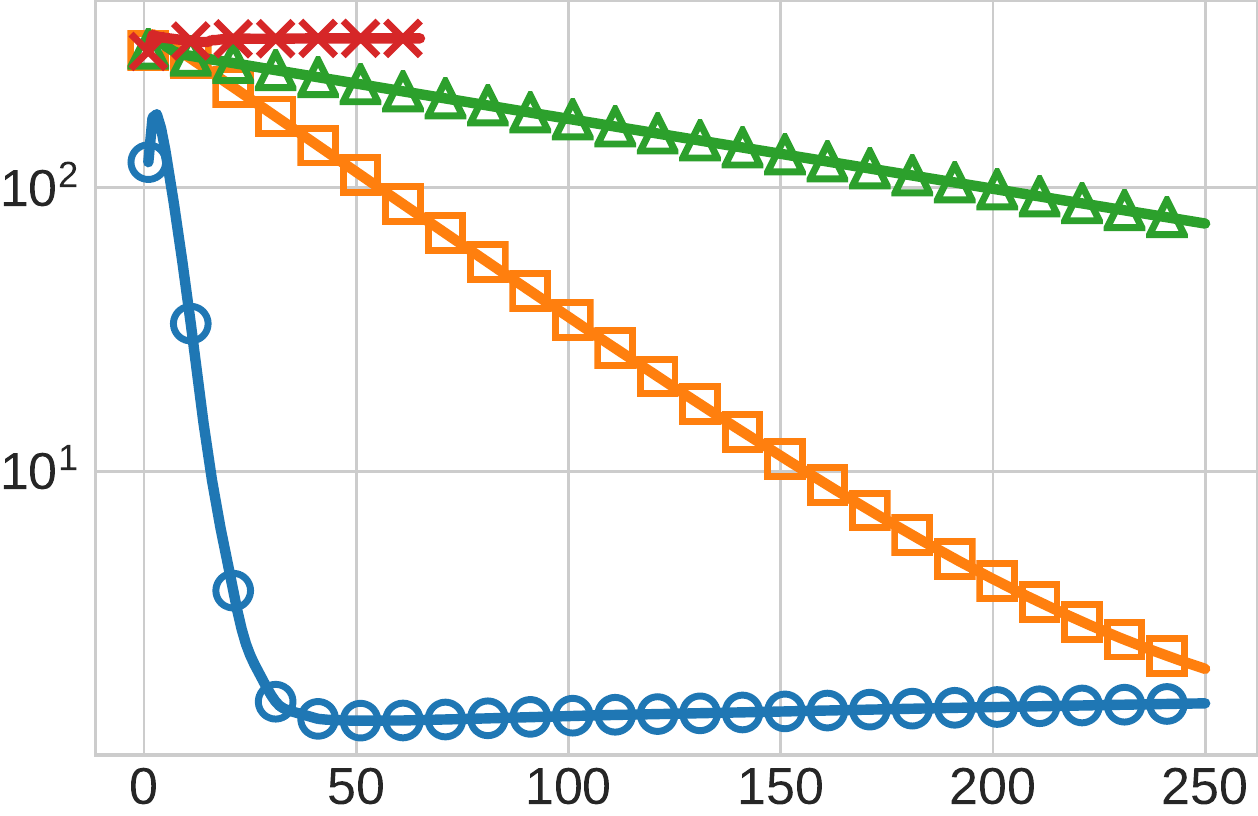}
			&
			\includegraphics[width=0.3\textwidth]{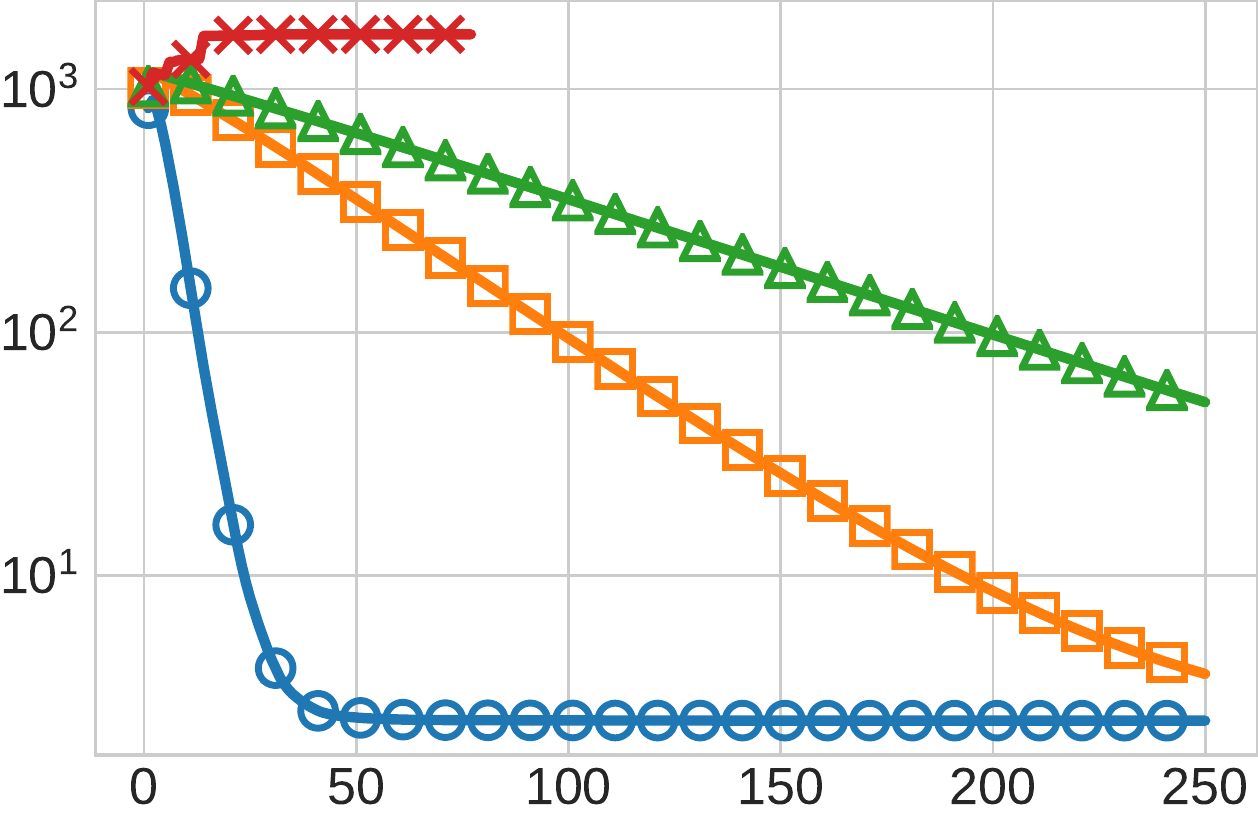} 
			&
			\includegraphics[width=0.3\textwidth]{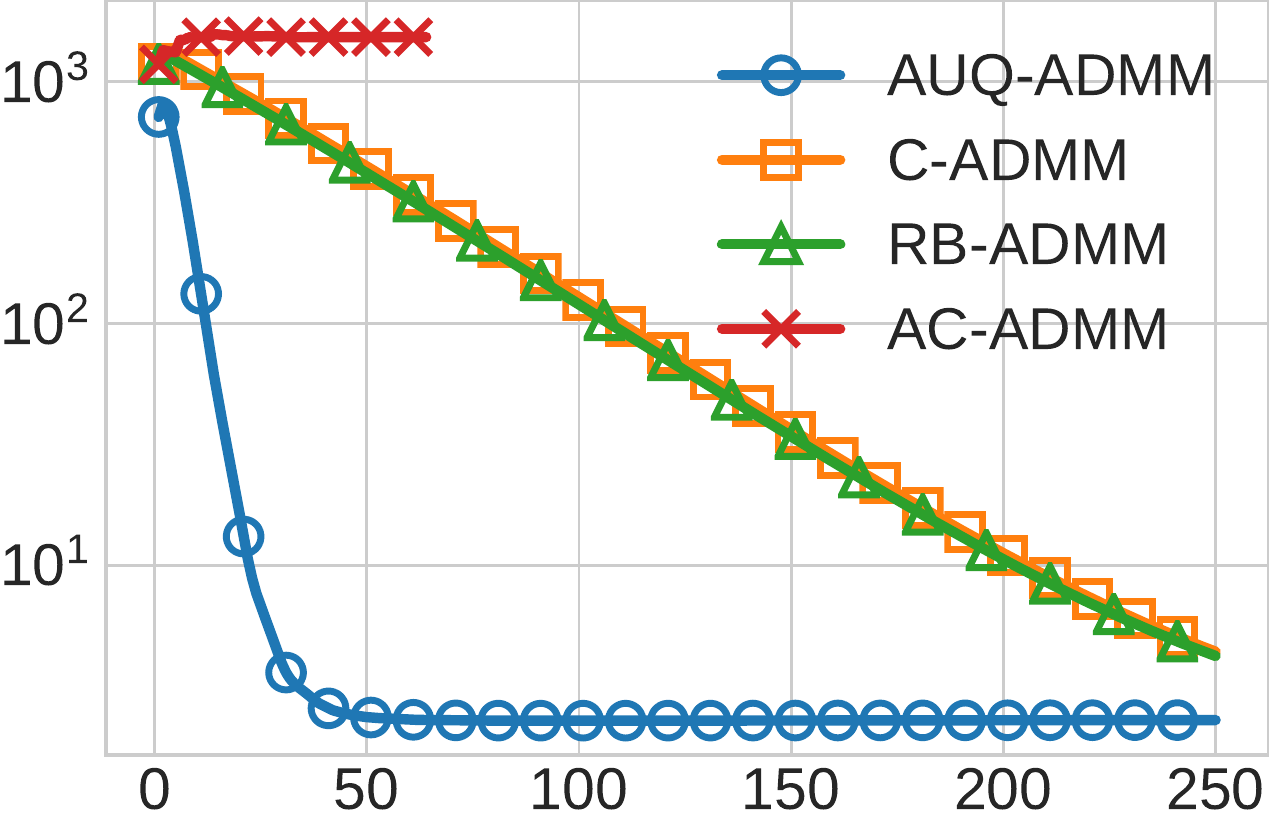}
			\\
			iterations & iterations & iterations
		\end{tabular}
        \caption{Multinomial Logistic Regression with MNIST, SVHN and CIFAR10: loss function value comparison among different ADMM-algorithms.}
		\label{fig:Multinomial}
	\end{figure}

    We consider the multinomial logistic regression with Tikhonov regularization. 
    In particular, $f_j(\bfu_j)$ and $g(\bfv)$ are defined as:
    \begin{equation}
        \label{eq: MultiReg}
        f_j(\bfu_j) = -\sum_{i=1}^{n_j}\sum_{k=1}^{C}\mathds{1}\{\bfy_{i}=k\}\log\frac{\exp{\left((\bfu_{j,k})^\top\bfX_{i}\right)}}{\sum_{l=1}^{C}\exp{\left((\bfu_{j,l})^\top\bfX_i\right)}},
        \quad g(\bfv)=\frac{1}{2}\|\bfv\|^2_2,
    \end{equation}
    where $\bfu_j,\bfv\in\bbR^{m\times C}$ in which $C$ is the number of classes, $n_j$ is the sample size for worker $j$, and $\bfu_{j,i}$ is the $i$th column of $\bfu_j$; $\mathbf{X}\in\bbR^{n_j\times m}$ is the feature matrix, $\mathbf{X}_i$ is the $i$th column of $\mathbf{X}$, $\bfy\in\bbR^{n_j}$ is the corresponding labels and $\bfy_i\in\{0,1,\dots,C-1\}$ is the $i$th component of $\bfy$. 
    In this case, $C=10$; $m=784$, $n_j=2000$ for MNIST; $m=3072$, $n_j=1000$ for SVHN and CIFAR10 as before. 
    The convergence results are shown in Figure \ref{fig:Multinomial}.
	
    
    \subsection{Support Vector Machines}
    \label{subsec: support_vector_machines}
\begin{figure}[htb]
		\centering
		\setlength{\tabcolsep}{1pt}
		\small
		\begin{tabular}{ccc}
		    MNIST Loss $f+g$ & SVHN Loss $f+g$ & CIFAR10 Loss $f+g$
			\\
			\includegraphics[width=0.3\textwidth]{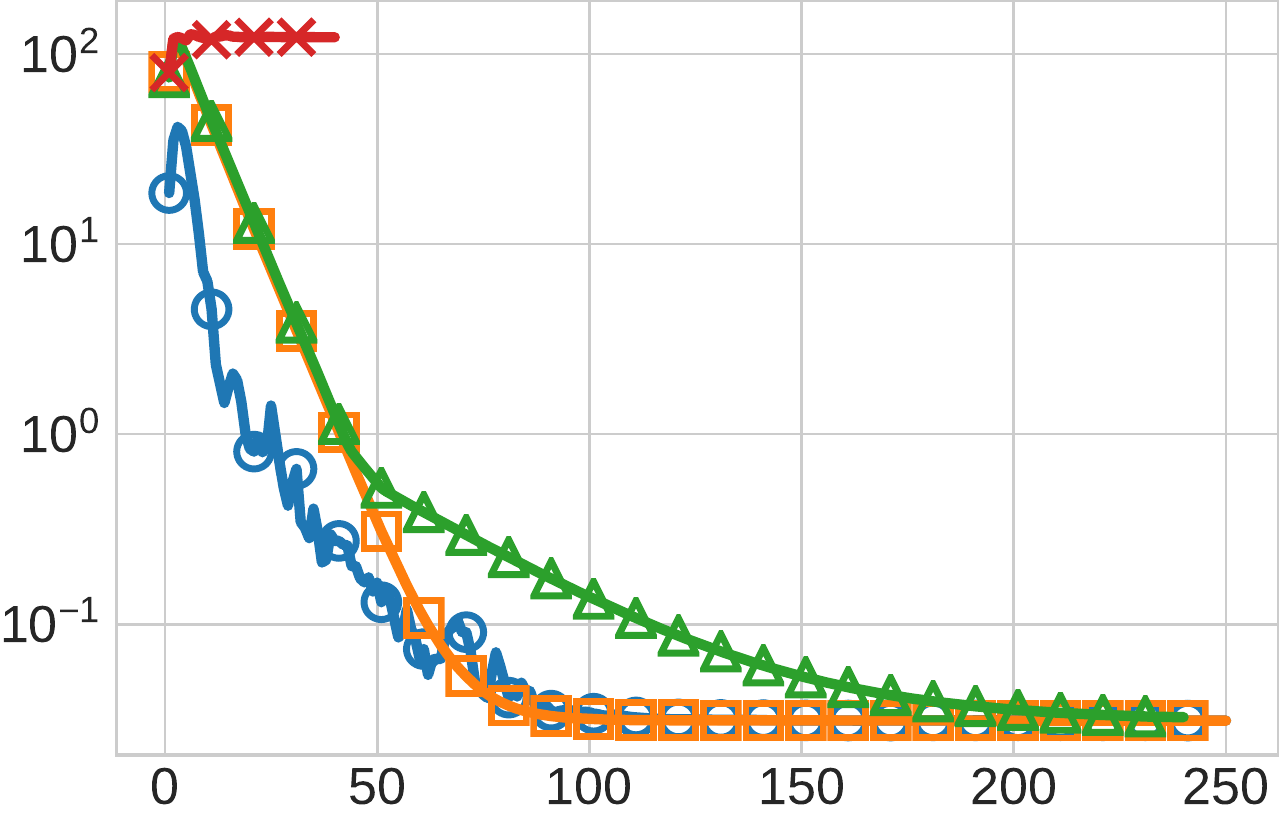}
			&
			\includegraphics[width=0.3\textwidth]{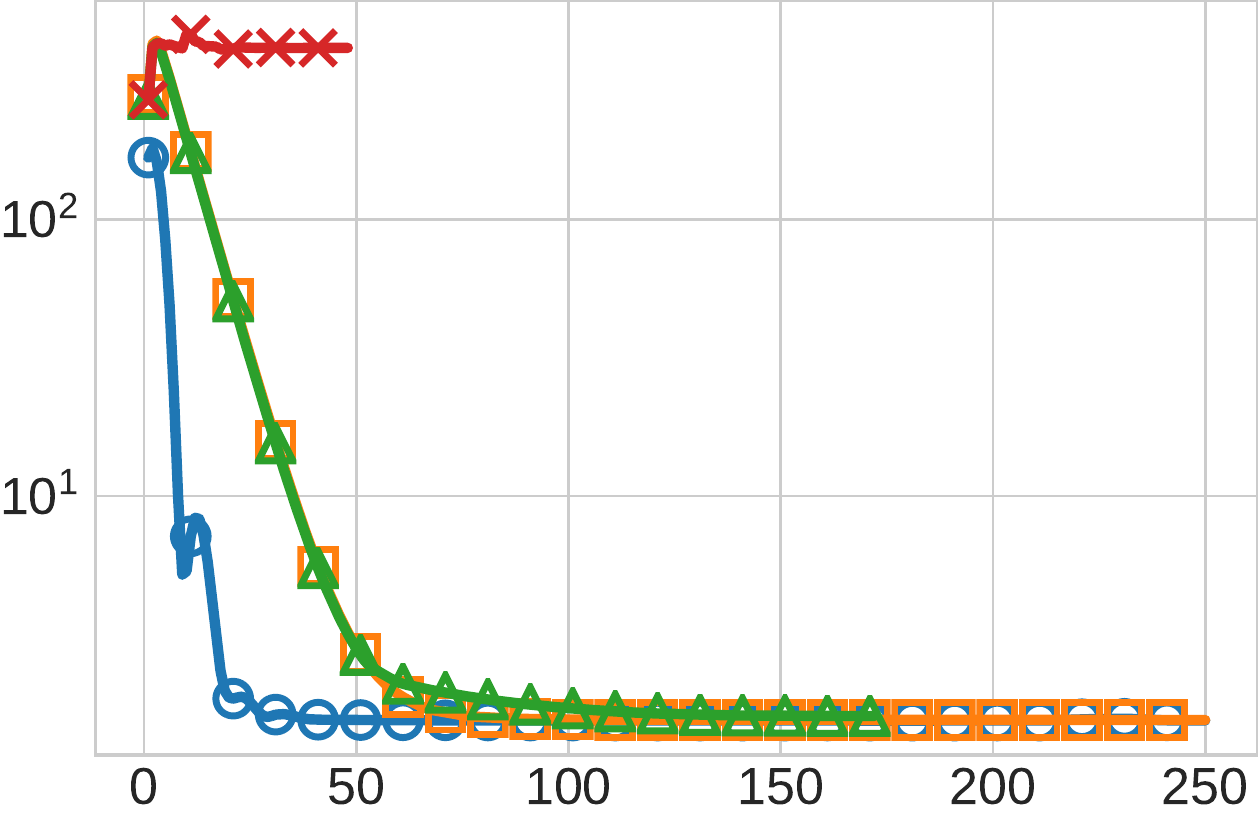} 
			&
			\includegraphics[width=0.3\textwidth]{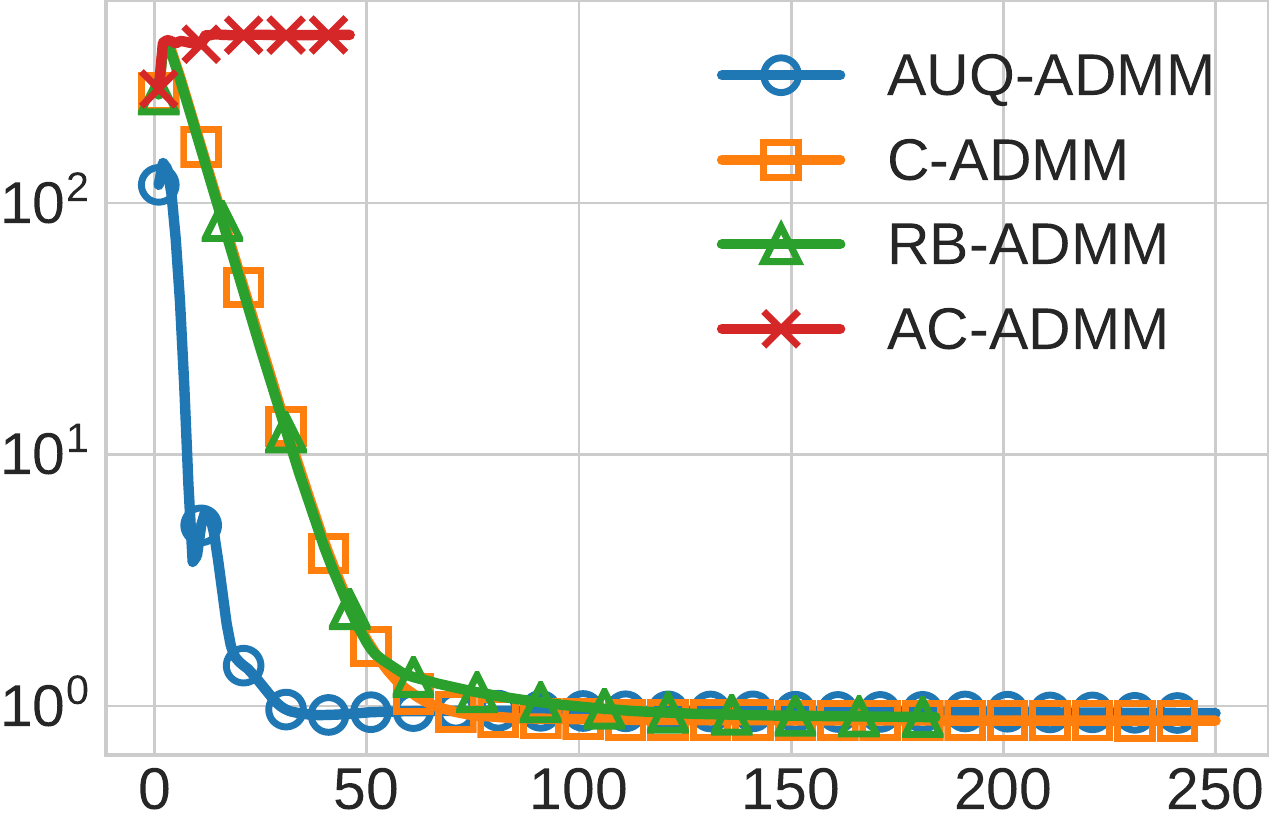}
			\\
			iterations & iterations & iterations
		\end{tabular}
        \caption{Smoothed Support Vector Machine with MNIST, SVHN and CIFAR10: loss function value comparison among different ADMM-algorithms.}
		\label{fig:SVM}
	\end{figure}

    We consider the soft-margin support vector machines (SVM)~\cite{8637389} for binary classification, with $f_j(\bfu_j)$ and $g(\bfv)$ defined as
    \begin{equation}
        \label{eqn: SVM}
        f_j(\bfu_j) = \frac{1}{n_j}\sum_{i=1}^{n_j}\frac{1}{2}\left( z_i+\sqrt{\epsilon^2+z_i^2} \right),
        \quad g(\bfv) = \frac{1}{2}\|\bfv\|^2.
    \end{equation}
    Here, 
    $$z_i = 1 - \bfy_i(\bfu_j)^{\top}\bfX_i,\quad \epsilon=\frac{1}{5000},$$
    and $\bfu_j,\bfv\in\bbR^{m}$, $n_j$ is the sample size for worker $j$; $\mathbf{X}\in\bbR^{m\times n_j}$ is the input data matrix, $\mathbf{X}_i$ is the $i$th column of $\mathbf{X}$, $\bfy\in\bbR^{n_j}$ is the corresponding labels and $\bfy_i\in\{-1,1\}$ is the $i$th component of $\bfy$. 
    In our case, $m=784$, $n_j=5000$ for MNIST; $m=3072$, $n_j=5000$ for CIFAR10; $m=3072$, $n_j=4500$ for SVHN. 
    Convergence histories are shown in Fig.~\ref{fig:SVM}.
    
    \subsection{Number of Splittings Comparison}
    A key feature of consensus ADMM is that performance of the method deteriorates as the number of splittings increases. 
    The intuition is that each worker has less information leading to slow convergence.
    In this experiment, we demonstrate the robustness of AUQ-ADMM when the number of splittings increase.
    We compare all algorithms using 2, 4, and 8 workers for the MNIST dataset. 
    For the convenience of splitting the dataset evenly, we only used 8 classes of MNIST (0 to 7) to perform this experiment. 
    \begin{figure}[htb]
		\centering
		\setlength{\tabcolsep}{1pt}
		\small
		\begin{tabular}{ccc}
		    {\normalsize 2 Workers} & {\normalsize 4 Workers} & {\normalsize 8 Workers}
			\\
			Loss $f+g$ & Loss $f+g$ & Loss $f+g$
			\\
			\includegraphics[width=0.3\textwidth]{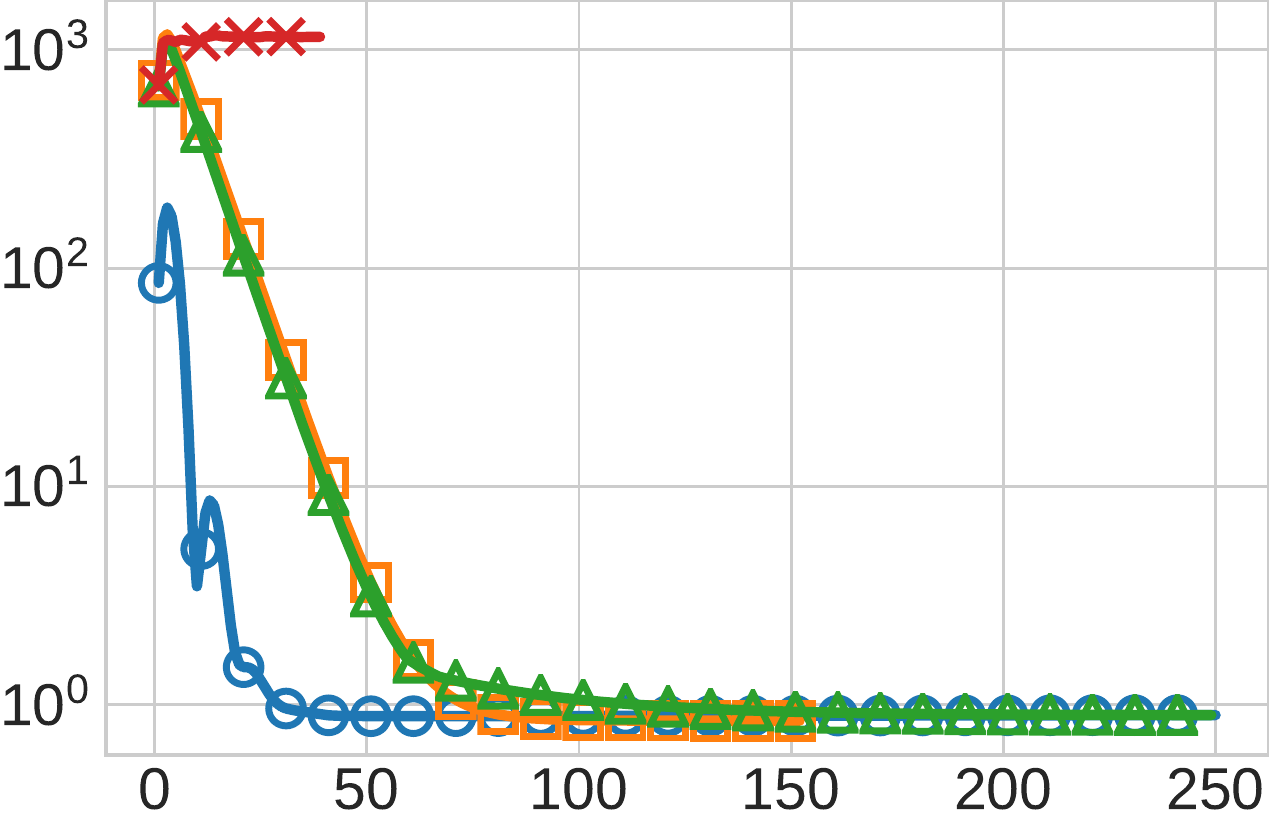}
			&
			\includegraphics[width=0.3\textwidth]{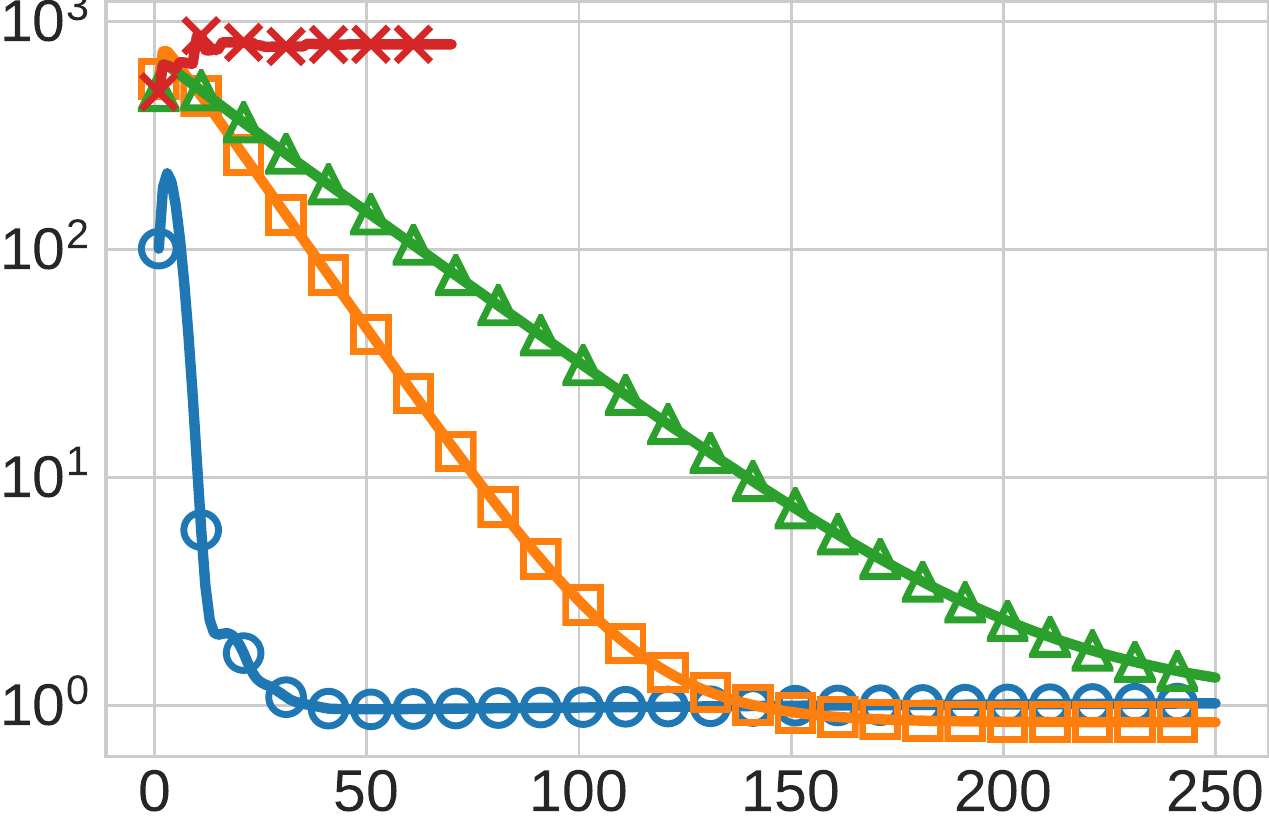} 
			&
			\includegraphics[width=0.3\textwidth]{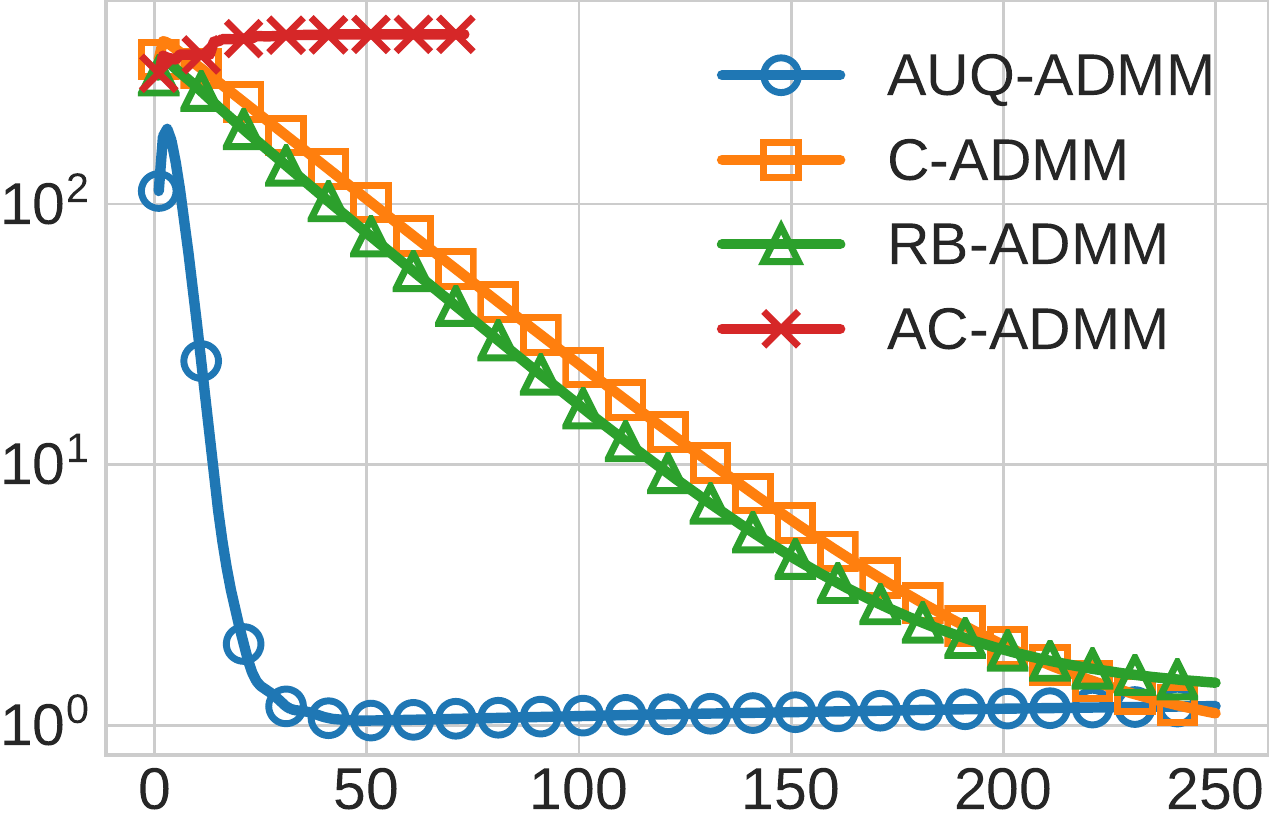}
			\\
			iterations & iterations & iterations
		\end{tabular}
		\caption{Convergence of AUQ-ADMM and C-ADMM under different number of workers, with multinomial logistic regression on MNIST dataset.}
		\label{fig:Workers}
	\end{figure}
	
	For 8 workers, each worker handles one class, with 2000 samples per worker.
	For 4 workers, each worker handles two classes, with 4000 samples per worker.
	For 2 workers, each worker handles 4 classes, 8000 samples per worker. 
	

	

\subsection{\edits{Uncertainty-based Weights}}
\label{subsec: weight_plots}
\begin{figure}[!htb]
    \centering
    \setlength{\tabcolsep}{0.2pt}
    \begin{tabular}{cccccccccc}
        $\bfW_1$ & $\bfW_2$ & $\bfW_3$ & $\bfW_4$ & $\bfW_5$ & $\bfW_6$ & $\bfW_7$ & $\bfW_8$ & $\bfW_9$ & $\bfW_{10}$
        \\
        \includegraphics[width=0.09\textwidth]{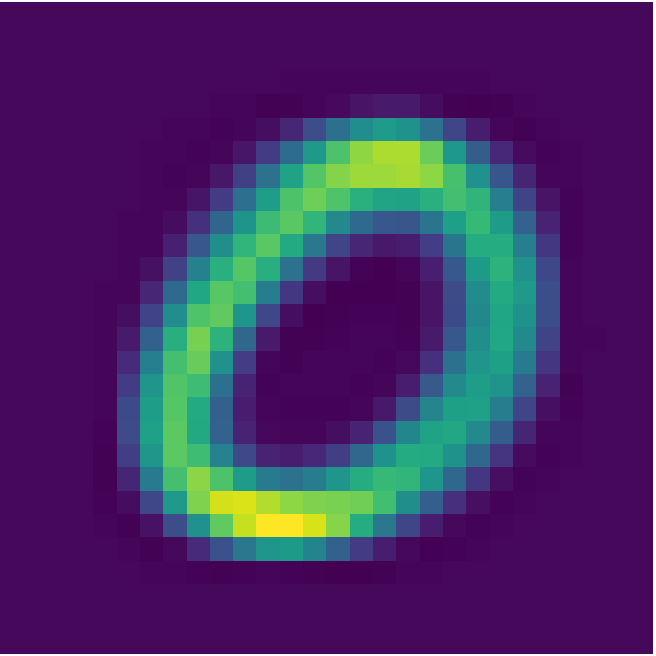}
        &
        \includegraphics[width=0.09\textwidth]{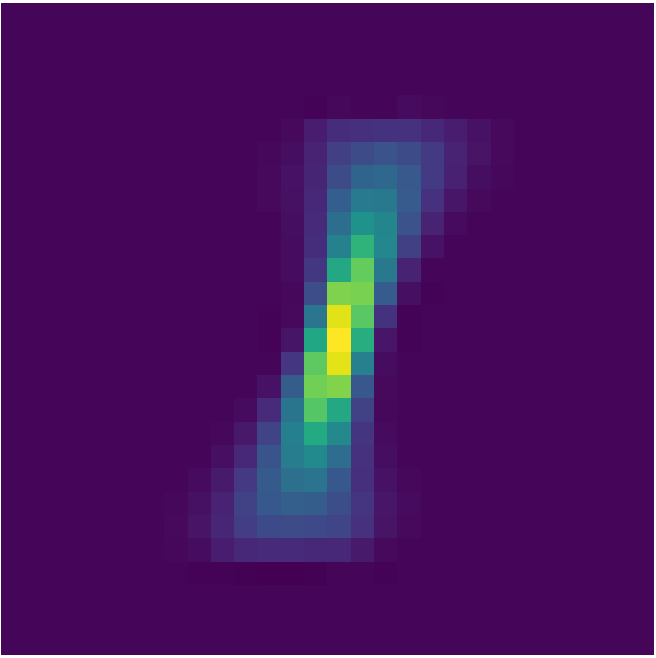}
        &
        \includegraphics[width=0.09\textwidth]{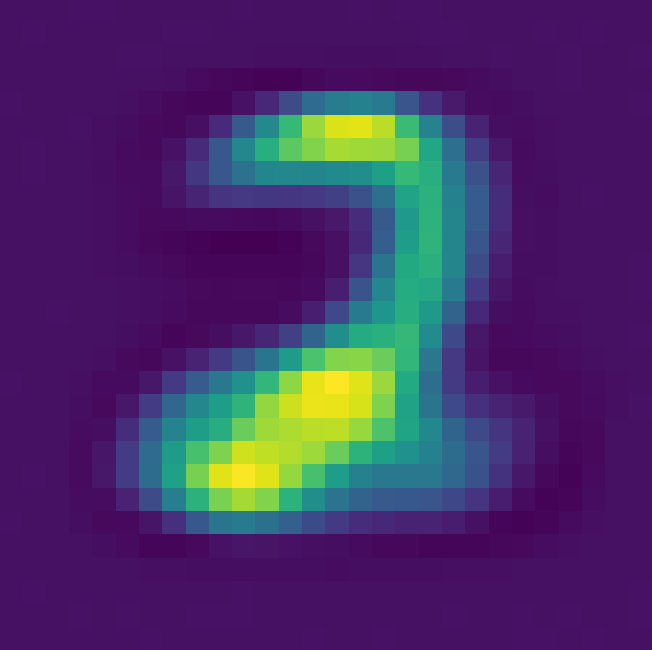}
        &
        \includegraphics[width=0.09\textwidth]{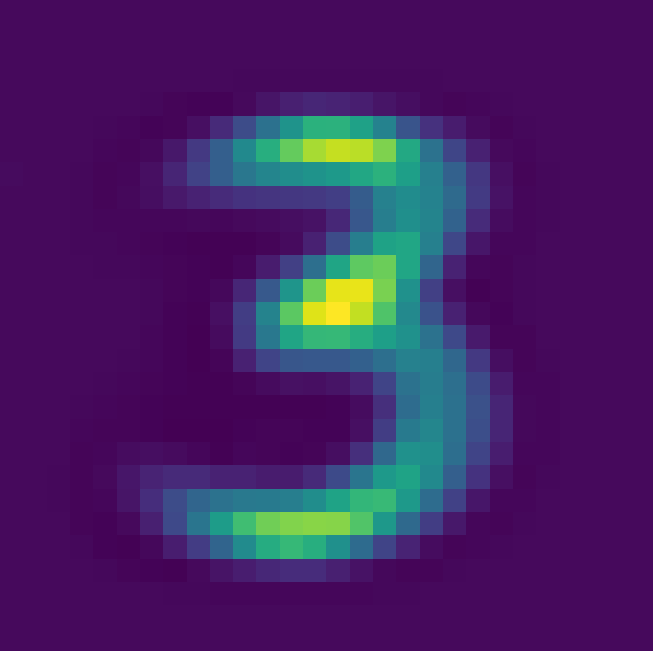}
        &
        \includegraphics[width=0.09\textwidth]{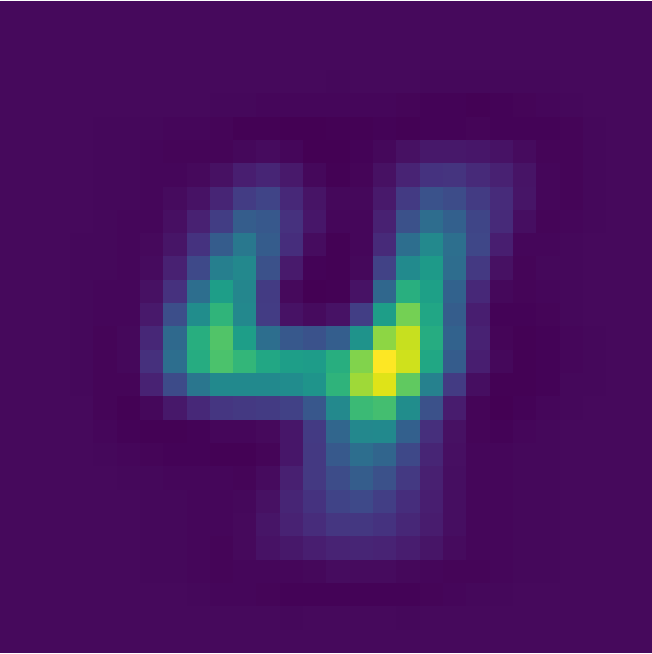}
        &
        \includegraphics[width=0.09\textwidth]{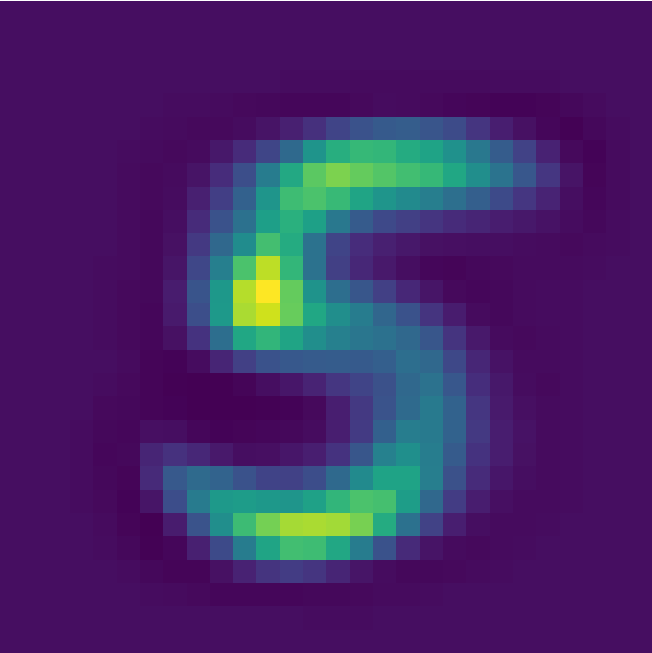}
        &
        \includegraphics[width=0.09\textwidth]{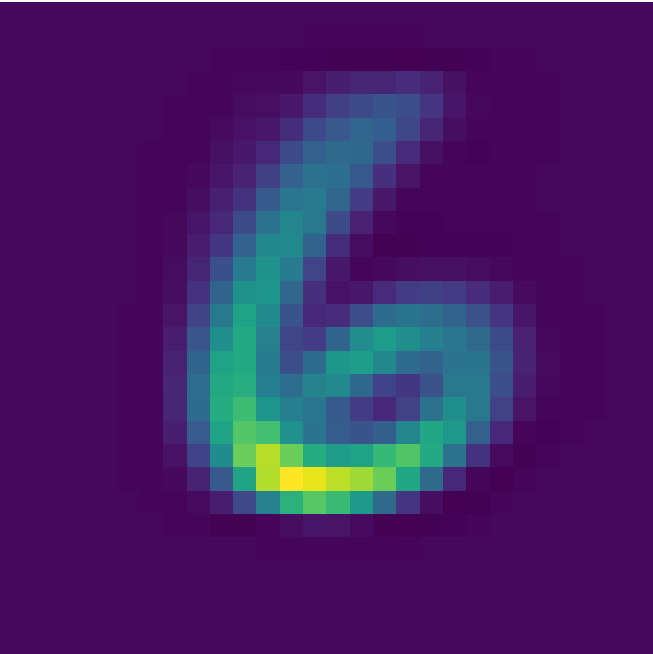}
        &
        \includegraphics[width=0.09\textwidth]{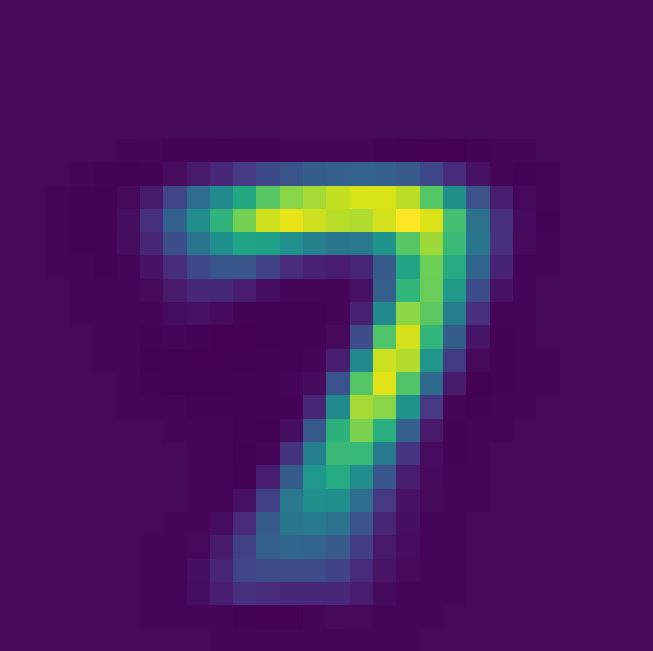}
        &
        \includegraphics[width=0.09\textwidth]{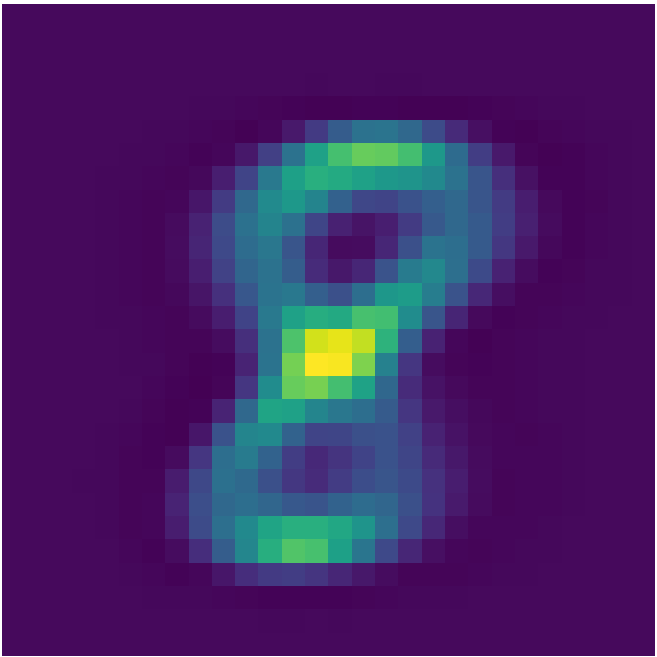}
        &
        \includegraphics[width=0.09\textwidth]{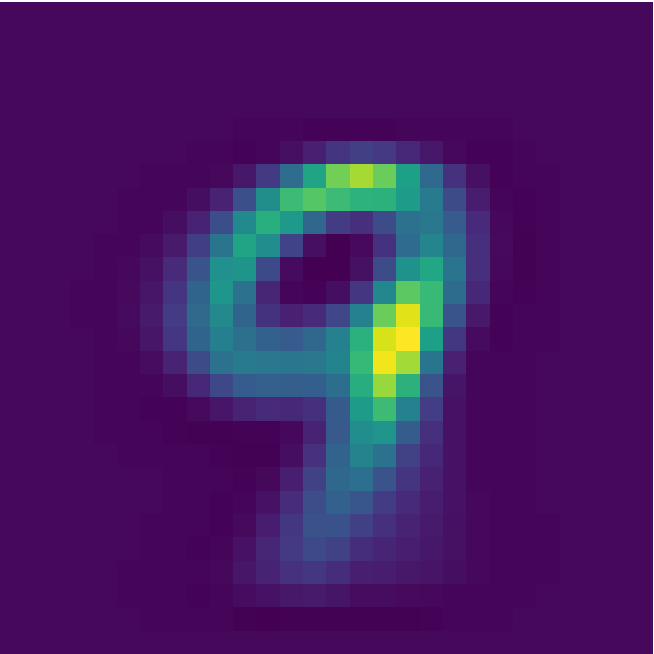}
    \end{tabular}
    \caption{Plot of \textit{diagonal} elements of weights for each class with elastic net loss.}
    \label{fig:weights_elasticnet}
\end{figure}

\begin{figure}[!htb]
    \centering
    \includegraphics[width=0.85\textwidth]{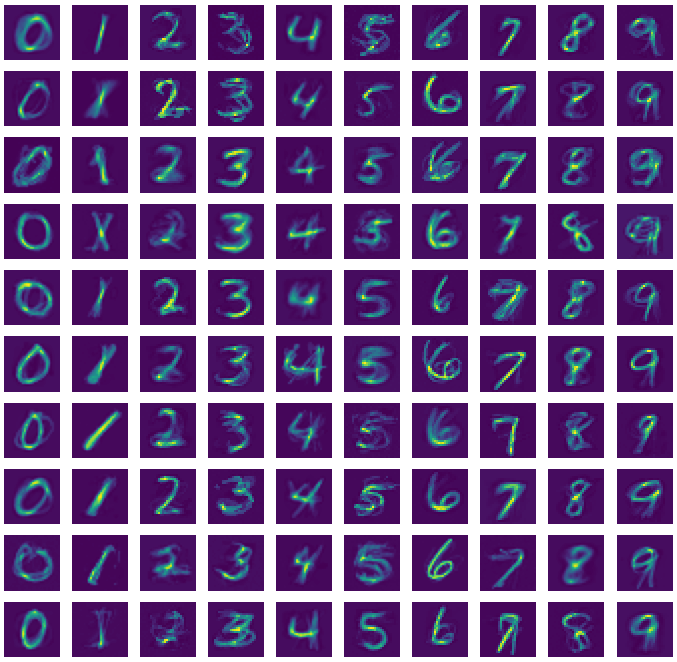}
    \caption{Plot of \textit{diagonal} elements of weights for each class with multinomial loss.}
    \label{fig:weights_multinomial}
\end{figure}

\begin{figure}[!htb]
    \centering
    \begin{tabular}{cc}
        $\bfW_1$ & $\bfW_2$
        \\
        \includegraphics[width=0.3\textwidth]{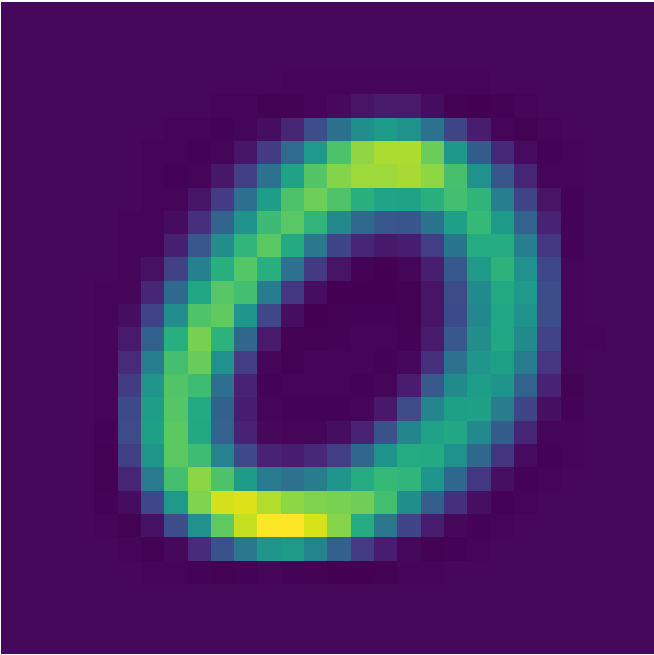}
        &
        \includegraphics[width=0.3\textwidth]{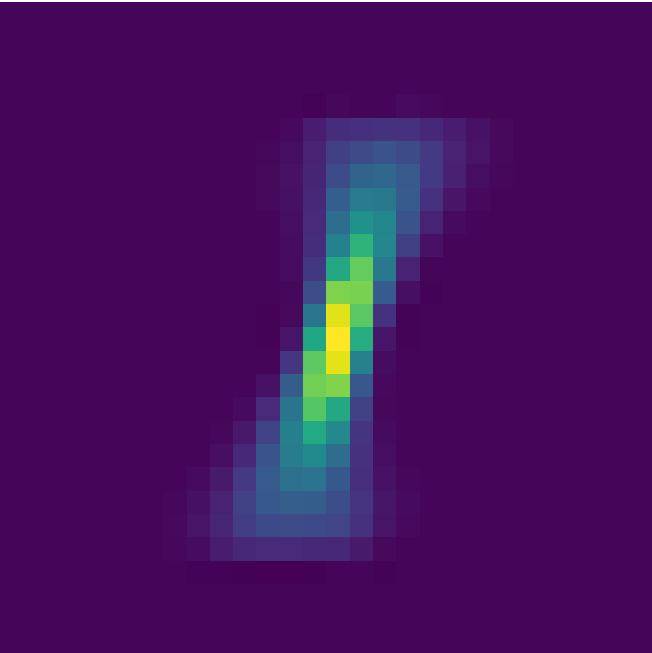}
    \end{tabular}
    \caption{Plot of \textit{diagonal} elements of weights for each class with SVM loss.}
    \label{fig:weights_svm}
\end{figure}

\edits{
To better understand the role of the uncertainty-based weights in our experiments, we plot the weights for the MNIST dataset for the different experiments presented in Sec.~\ref{subsec: elastic_net_regression},~\ref{subsec: multinomial_logistic_regression}, and~\ref{subsec: support_vector_machines}.
Fig.~\ref{fig:weights_elasticnet} shows the weights when $f_j$ is based on elastic net regression. In particular, we see that $\bfW_1$ resembles an averaged image of all the zeros, $\bfW_2$ resembles an averaged image of all the ones, and so on.
This is aligned with our intuition as we do not wish our model $\bfu_j$ to be weighted where the pixels are zero-valued when updating the global variable $\mathbf{v}$. 
Similarly, we observe similar features in the weights for multinomial logistic regression in Fig.~\ref{fig:weights_multinomial}, and support vector machines in~\ref{fig:weights_svm}.
Note that the weights for Fig.~\ref{fig:weights_multinomial} contain ten images per class by design (as their size are $\bfW_j \in \bbR^{784 \times 10}$) - this is the reason why Fig~\ref{fig:weights_multinomial} contains 100 images, each of size $28 \times 28$.
Similarly, our SVM experiment was performed on a binary classification problem, and thus, there are only two images.
Ultimately, we observe that the weights have low values in uninformative areas of the features and vice-versa, leading to improved averaged reconstructions, and ultimately, improved convergence.
}

\subsection{Discussion}
\label{subsec: discussion}
    For all experiments, AUQ-ADMM outperforms all other methods in decreasing the objective function Fig.~\ref{fig:Elastic}, Fig.~\ref{fig:Multinomial}, and Fig.~\ref{fig:SVM}. We also experimentally observe improved robustness to number of splittings in Fig.~\ref{fig:Workers}, which motivates using AUQ-ADMM for large datasets where many splittings are required.
    One reason for this the ``informed'' averaging performed during the global variable update as was illustrated in Fig.~\ref{fig:weightIllustration}. This averaging assigns higher weights to elements in the local models where we are more certain and vice-versa. \edits{The weights in Fig.~\ref{fig:weights_elasticnet},~\ref{fig:weights_multinomial}, and~\ref{fig:weights_svm} also support our intuition.}
    
\section{Conclusion}
    \label{sec: conclusion}
    We present AUQ-ADMM, an adaptive, weighted, consensus ADMM method for solving large-scale distributed optimization problems. 
    Our proposed weighting scheme is based on the uncertainties of the model.
    Intuitively, the uncertainty framework assigns higher weights where we are more certain and vice versa.
    The weights are computed efficiently using a low-rank approximation.
    Convergence of AUQ-ADMM is provided.
    
    Our experiments show that uncertainty-based weights improves performance for elastic net regression, multinomial logistic regression, and support vector machines using the MNIST, SVHN, and CIFAR10 datasets.
    Moreover, AUQ-ADMM is more robust to increased splittings. A natural extension of UQ-ADMM involves its application to more general non-convex problems such as training deep neural networks for classification~\cite{lecun2015deep}. In particular, the robustness with respect to the number of splittings could be well-tailored toward training swarm-based multi-agent control models~\cite{onken2021neural,onken2021neural_full}.
    However, convergence guarantees for the non-convex setting is a more difficult task and is a direction we intend to pursue.
%
%

\vskip 6mm
\noindent{\bf Acknowledgments}

\noindent We thank Lars Ruthotto for the fruitful discussions. We also thank the anonymous referees for useful suggestions which improved the contents of this paper.
\bibliographystyle{spmpsci}      
\bibliography{references}   

\newpage
\appendix

\section{Convergence Proofs}
\label{Appendix: A}
This section provides proofs for the results of Sec.~\ref{sec: convergence}. For the reader's convenience, we restate the lemmas and theorems before proving them. Once Lemma~\ref{lemma: WeightConvg} is proven, the remaining results are a straightforward application of the theorem proven in~\cite{xu2017adaptive2}.
~\\
{\bf Lemma \ref{lemma: WeightConvg}.} \textit{
In Algorithm~\ref{alg:auqdmm}, the weights $\{\bfW^{(k)}\}_{k=0}^{\infty}$ ($\bfW^{(0)}=\bfW^{(1)}$) generated by AUQ-ADMM satisfy
        \begin{align}
            \bfW^{(k)} \preceq (1+c^{(k)}) \bfW^{(k-1)} \; \text{ and } \; (\bfW^{(k)})^{-1} \preceq (1+c^{(k)}) (\bfW^{(k-1)})^{-1}, 
          \end{align}
        where
        \begin{equation}
        c^{(k)} = 
        \begin{cases}
          \frac{b_1}{a_1}-1& \text{if $k=1$}
          \\
          \frac{b_{k-1}}{a_{k-1}}-1. & \text{if $k>1$}
        \end{cases},
        \qquad k = 1,2,3,\ldots
        \end{equation}
        is a sequence of positive scalars satisfying
        \begin{equation}
            \sum_{k=1}^\infty c^{(k)} < \infty \notag,
        \end{equation}
        and $[a_1, b_1]$ is the initial restriction interval. 
}
\paragraph{\textbf{Proof of Lemma \ref{lemma: WeightConvg}}:}
Denote the initial interval given by the user as $[a_1, b_1]$ and the $i$th updated interval as $[a_i,b_i]$. We define the following sequence $c^{(k)}$,
\begin{equation}
        \label{eqn: ck_appdix}
          c^{(k)} =
          \begin{cases}
              \frac{b_1}{a_1}-1& \text{if $k=1$}\\
              \frac{b_{k-1}}{a_{k-1}}-1. & \text{if $k>1$}
          \end{cases}
\end{equation}

We first show $\sum_{k=1}^{\infty}c^{(k)}<\infty$. From Algorithm \ref{alg:IntUpdt} we can observe that for every $i=1,2,3,\dots$,
\begin{equation}
    \label{important_ineqty}
    a_i=a_1,\quad \frac{b_i}{a_i}=\frac{1}{i^2}\frac{b_1}{a_1}+1-\frac{1}{i^2},
\end{equation}

then,
\begin{align}
    \sum_{k=1}^{\infty}c^{(k)} &=\frac{b_1}{a_1}-1+ \sum_{i=1}^{\infty} \left(\frac{b_i}{a_i}-1\right) \notag \\
    &= \frac{b_1}{a_1}-1+\sum_{i=1}^{\infty}\left(\frac{1}{i^2}\frac{b_1}{a_1}-\frac{1}{i^2}\right)<\infty. \notag
\end{align}
Now we show inequality (\ref{tempEQ}) holds with this $c^{(k)}$. For convenience, let $w_j^{(k)}$ denote the $j^{\textit{th}}$ diagonal element of $\bfW^{(k)}$. Observe that showing $(\ref{tempEQ})$ is equivalent to showing
\begin{equation}
\label{inequality}
    \frac{1}{1+c^{(k)}}w^{(k-1)}_j\leq w^{(k)}_j \leq(1+c^{(k)})w^{(k-1)}_j,\quad \text{for every }j=1,2,\dots,Nn.
\end{equation}

And we define the following notations:
\begin{equation}
    q_{\max} = \max_{1\leq j\leq Nn}\frac{w^{(k-1)}_j}{w^{(k)}_j},\quad q_{\min} = \min_{1\leq j\leq Nn}\frac{w^{(k-1)}_j}{w^{(k)}_j}. \notag
\end{equation}

As $\bfW^{(0)}=\bfW^{(1)}$, it is obvious that (\ref{inequality}) holds when $k=1$. When $k>1$, the restriction interval update algorithm tells us that we have restriction interval $[a_k,b_k]$ at $k$th step, but $[a_{k+1},b_{k+1}]$ at $(k+1)$th step. Also, we observe that,
$$q_{\max}\leq\frac{b_{k}}{a_{k+1}}=\frac{b_k}{a_k}, \quad q_{\min}\geq\frac{a_{k}}{b_{k+1}}\geq\frac{a_k}{b_k},$$
because $a_{k+1}=a_k=a_1$ and $b_{k+1}\leq b_k$ (this fact can be derived from (\ref{important_ineqty})). 

Since $c^{(k+1)}=\frac{b_k}{a_k}-1$, we can derive that
$$(1+c^{(k+1)})^2=\frac{b_k^2}{a_k^2}\geq\frac{q_{\max}}{q_{\min}}\quad\Rightarrow\quad\frac{1}{1+c^{(k+1)}}q_{\max}\leq(1+c^{(k+1)})q_{\min}.$$

In addition, note that
$$(1+c^{(k+1)})q_{\min}\geq\frac{b_k}{a_k}\frac{a_k}{b_k}=1$$
and
$$\frac{1}{(1+c^{(k+1)})}q_{\max}\leq\frac{a_k}{b_k}\frac{b_k}{a_k}=1.$$

Therefore, we can still have
$$\frac{1}{1+c^{(k+1)}}\frac{w^{(k)}_j}{w^{(k+1)}_j}\leq\frac{1}{1+c^{(k+1)}}q_{\max}\leq1\leq(1+c^{(k+1)})q_{\min}\leq(1+c^{(k+1)})\frac{w^{(k)}_j}{w^{(k+1)}_j}.$$

Combining the results completes the proof of Lemma \ref{lemma: WeightConvg}.
$\hfill \square$
\\

\begin{lemma}
\label{essentiallem}
For any $\bfy=(\bfu;\bfv;\bflambda)$ and $\bfy'=(\bfu';\bfv';\bflambda')$, we have,
\begin{equation}
    \|\bfy-\bfy'\|^2_{\bfT^{(k)}}\leq (1+c^{(k)})\|\bfy-\bfy'\|^2_{\bfT^{(k-1)}} \notag
\end{equation}
\end{lemma}
\textit{Proof:} This immediately follows from (\ref{tempEQ}). $\hfill \square$
\\

The following theorem immediately follows from Theorem 1 in AC-ADMM paper \cite{xu2017adaptive2}:
\begin{theorem}
\label{Thm1}
The sequence $\bfy^{(k)}=(\bfu^{(k)};\bfv^{(k)};\bflambda^{(k)})$ generated by the AUQ-ADMM satisfies
\begin{equation}
    \lim_{k\gets\infty}\|\bfy^{(k+1)}-\bfy^{(k)}\|^2_{\bfT^{(k)}}=0.
\end{equation}
\end{theorem}

Theorem \ref{theorem: Convg} in Sec.~\ref{sec: convergence} is just a simple variant of Theorem 2 in the AC-ADMM paper by replacing $z$ with optimal $z^*$.

Here are the interpretations of these theorems: Theorem \ref{Thm1} shows the convergence of the scheme satisfying Assumption 1, and Theorem \ref{theorem: Convg} gives a specific error bound for the ergodic $O(1/k)$ convergence rate.

\end{document}